
\input amstex
\documentstyle{amsppt}
\document
\magnification=1200
\NoBlackBoxes
\nologo
\hoffset=0.7in
\voffset=1in
\pageheight{16cm}


\bigskip

{\bf
\centerline{QUANTUM THETA FUNCTIONS} 
\smallskip
\centerline {AND}
\smallskip
\centerline{ GABOR FRAMES FOR MODULATION SPACES}}

\bigskip

\centerline{\bf Franz Luef${}^{(1)}$, Yuri I. Manin${}^{(2)}$}

\medskip

\centerline{\it ${}^{(1)}$ University of Vienna, Vienna, Austria}
\centerline{\it and University of California, Berkeley, USA}
\smallskip

\centerline{\it ${}^{(2)}$ Max--Planck--Institut f\"ur Mathematik, Bonn, Germany}
\centerline{\it and Northwestern University, Evanston, USA}

\bigskip

{\bf Abstract.} Representations of the celebrated
Heisenberg commutation relations in quantum mechanics (and their
exponentiated versions) form the starting point
for a number of basic constructions, both in mathematics and mathematical physics
(geometric quantization, quantum tori, classical and quantum theta functions) 
and signal analysis (Gabor analysis). In this  paper we will try to bridge the two communities,
represented by the two co--authors: that of noncommutative
geometry  and that of signal analysis. After providing a brief comparative dictionary of the two languages,
we will show e.g. that the Janssen representation of Gabor frames with generalized Gaussians as Gabor atoms 
yields in a natural way quantum theta functions, and that the Rieffel
scalar product and associativity relations underlie both the functional equations for quantum thetas
and the Fundamental Identity of Gabor analysis.

\bigskip

\centerline{\bf \S 0. Introduction.}

\medskip

Gabor analysis is a modern branch of signal analysis with various 
applications to pseudodifferential operators, harmonic analysis, function spaces, approximation theory, and quantum mechanics. 
It is  well known  that there are substantial connections between the mathematical foundations of signal analysis and those of 
quantum mechanics. Furthermore, the theory of operator algebras furnished  a  rigorous framework for quantum mechanics, 
but possible direct relationships between signal analysis and operator algebras have not received much attention.
Recent work of Gr\"ochenig and his 
collaborators (see [Gr], [GrLe]) made explicit connections between the spectral invariance of certain Banach algebras and basic problems in signal analysis. 
In his Ph.D. thesis [Lu1] one of us has developed basic correspondences between Gabor analysis and noncommutative geometry over noncommutative tori. 

\smallskip

Classical theta functions have a long history, see [Mu] for a modern exposition.   
From the functional theoretic viewpoint, they are
holomorphic functions of several complex variables, which acquire an exponential factor
after the shift by any vector in a period lattice. Geometrically, they represent
sections of a line bundle over a complex torus, lifted to the universal cover of this torus
and appropriately trivialized there.

\smallskip

The second named author, motivated by the ideas of geometric quantization, suggested 
in 1990 that one can develop a meaningful theory of quantum thetas
after replacing ordinary complex tori in the classical construction by their
quantum versions, see [Ma1]. The new quantum theta functions were subsequently applied to 
the construction of algebraic quantization of abelian varieties in [Ma2]
(this case, as well as that of symplectic projective manifolds
in general, presented a problem in Kontsevich's  paper [Ko]) and to the program 
of Real Multiplication ([Ma4], [Vl]). Thanks to the Boca study [Bo1],
it became clear that quantum thetas can be constructed as Rieffel's scalar products of
vacuum vectors in representations of the appropriate Heisenberg groups.
This idea, developed in [Ma5], led to the discovery 
of a quantum version of those classical functional equations
for theta functions that arise from different natural trivializations of
a line bundle over a complex torus, see also [EeK].

\smallskip

In this paper  we  survey  a new interpretation of quantum theta functions in the framework of Gabor analysis. Recent  investigations have clarified and enriched parts of both subjects, see  [Lu4], [Lu5]
 for work related to Gabor frames, and [GrLu] 
for a contribution on the structure of projective modules over noncommutative tori relying on methods from Gabor analysis. The present work is another 
instance for the relevance of Gabor analysis in exploiting basic notions of noncommutative geometry. The basic link between Gabor analysis and noncommutative geometry is furnished by the Heisenberg group, especially the Schr\"odinger representation of the Heisenberg group. The Heisenberg group lies at the heart of various branches of physics, 
applied and pure mathematics, see the excellent survey [Ho]. After the groundbreaking work of A.~Weil [We], theta functions have been linked with 
representations of Heisenberg groups. In this famous paper  Weil introduced the {\it metaplectic representation}, which had independently been found by Shale. 
Weil's  new methods and objects have influenced many mathematicians in their work on theta functions, most notably Cartier, Igusa and Reiter in  [Ca],  [Ig], [Re1]--[Re3].
In his work on abelian varieties Mumford had demonstrated the relevance of the Heisenberg group in the
algebraization  of theta functions, cf. [Mu]. In [Sc],  Schempp has discussed the close relation between signal analysis and 
theta functions, where the Heisenberg group and its representation theory serves as a link between these two objects. The present investigation might be considered as a far--reaching extension of this line of research.

\bigskip

\centerline{\bf \S 1. Central extensions and Heisenberg groups.}

\medskip
In this section we recall the basic definitions of central extensions and Heisenberg groups. Our presentation 
follows closely the one given in \cite{Ma5}.
\smallskip  

{\bf 1.1. Central extensions.} Let $\Cal{K}$ (resp. $\Cal{Z}$)
be an abelian group written additively (resp. multiplicatively).
Consider a function $\psi :\,\Cal{K}\times \Cal{K}\to \Cal{Z}$.
Then the following conditions (a) and (b) are equivalent:

\smallskip

(a) {\it $\psi (0,0)=1$ and $\psi$ is a cocycle, that is, for each $x,y,z$
we have}
$$
\psi (x,y)\psi (x+y,z) =\psi (x,y+z)\psi (y,z) .
\eqno(1.1)
$$

(b) {\it The following composition law on $\Cal{G}:=\Cal{Z}\times\Cal{K}$
turns $\Cal{G}=\Cal{G} (\Cal{K},\psi)$ into a group with identity $(1,0)$:}
$$
(\lambda ,x)(\mu ,y):= (\lambda\mu\psi (x,y), x+y).
\eqno(1.2)
$$

Moreover, if (a), (b) are satisfied, then the maps
$\Cal{Z}\to\Cal{G}:\,\lambda\mapsto (\lambda ,0),$
$\Cal{G}\to \Cal{K}:\, (\lambda ,x)\mapsto x$, 
describe $\Cal{G}$ as a {\it central extension} of
$\Cal{K}$ by $\Cal{Z}$:
$$
1\to \Cal{Z}\to \Cal{G}(\Cal{K},\psi )\to \Cal{K}\to 1.
\eqno(1.3)
$$
Notice that any bicharacter $\psi$ automatically satisfies (a).
For arbitrary $\psi$, putting $x=0$ in (1.1), we see that
$\psi (0,y)=1$ so that 
$$
(\lambda ,x) =(\lambda ,0)(1,x).
$$

\smallskip

{\bf 1.1.1. Bicharacter $\varepsilon$.} Conversely, consider any central extension
(1.3), choose a set theoretic section $\Cal{K}\to \Cal{G}:\,x\mapsto
\tilde{x}$ and define the map $\varepsilon :\,
\Cal{K}\times \Cal{K}\to \Cal{Z}$ by
$$
\varepsilon (x,y):= \tilde{x}\tilde{y}\tilde{x}^{-1}\tilde{y}^{-1}.
\eqno(1.4)
$$
Then $\varepsilon$ is a bicharacter which does not depend on
the choice of a section and which is antisymmetric: $\varepsilon (y,x)=
\varepsilon (x,y)^{-1}$, $\varepsilon(x,x) =1.$ In particular,
if $K\subset \Cal{K}$ is a subgroup liftable to $\Cal{G}$,
then $K$ is $\varepsilon$--isotropic.
\smallskip

For the group
$\Cal{G} (\Cal{K}, \psi )$, choosing $\tilde {x}= (1,x)$,
we find
$$
\varepsilon (x,y) =\frac{\psi (x,y)}{\psi (y,x)},
\eqno(1.5)
$$
and if $\psi$ itself is an antisymmetric bicharacter, then
$\varepsilon (x,y) =\psi (x,y)^2.$ 

\medskip

{\bf 1.1.2. Cohomological interpretation.} The class of $\psi$ in $H^2(\Cal{K},\Cal{Z})$ determines $\Cal{G}$ up to an isomorphism identical on $\Cal{K},\Cal{Z}$. 
This extension is abelian iff $\varepsilon$  is trivial in which case the
extensions are classified by elements of $\roman{Ext}^1(\Cal{K},\Cal{Z})$. The map $\psi\mapsto\varepsilon$ coincides with the second arrow in
the universal coefficients exact sequence
$$
\roman{Ext}^1(\Cal{K},\Cal{Z}) \to H^2(\Cal{K},\Cal{Z})\to \roman{Hom}~(\Lambda^2\Cal{K},\Cal{Z}).
$$
\medskip

{\bf 1.2. Representations of central extensions.} Given $\Cal{K},\,\Cal{Z},\, \psi$
and a ground field $k$, choose in addition a character
$\chi :\, \Cal{Z} \to k^*$. Consider a linear space of functions
$\bold{f}:\,\Cal{K}\to k$ invariant with respect to the affine shifts 
and define operators $U_{(\lambda ,x)}$
on this space by
$$
(U_{(\lambda ,x)}\bold{f}) (x):= \chi (\lambda \psi(x,y)) \bold{f}(x+y).
\eqno(1.6)
$$ 
A straightforward check shows that this is a representation of
$\Cal{G} (\Cal{K}, \psi )$. However, it is generally reducible.
Namely, suppose that there is an $\varepsilon$--isotropic subgroup
$K_0 \subset \Cal{K}$ liftable to $\Cal{G} (\Cal{K}, \psi ).$ Let
$\sigma :\, K_0\to \Cal{G} (\Cal{K}, \psi ),$ $\sigma (y) = (\gamma (y), y)$
be such a lift. Denote by $F(\Cal{K}//K_0)$ the subspace
of functions satisfying the following condition:
$$
\forall\, x\in \Cal{K},\, y\in K_0,\
(U_{(\gamma (y) ,y)}\bold{f}) (x):= \chi (\varepsilon (x,y)) \bold{f}(x),
\eqno(1.7)
$$  
or, equivalently,
$$
\forall\, x\in \Cal{K},\, y\in K_0,\
\bold{f}(x+y) =\chi (\gamma (y)^{-1} \psi (y,x)^{-1}) \bold{f}(x).
\eqno (1.8)
$$
This subspace is invariant with respect to (1.6). If we choose for $K_0$ a maximal isotropic subgroup, then this provides 
a minimal subspace of this kind.

\smallskip

Formula (1.8) shows that if we know the value
of $\bold{f}$ at a point $x_0$ of $\Cal{K}$, it extends
uniquely to the whole coset $x_0+K_0$, hence the notation
$F(\Cal{K}//K_0)$ suggesting ``twisted'' functions
on the coset space $\Cal{K}/K_0$.

\medskip

{\bf 1.3. Locally compact topological groups.} The formalism briefly explained above is only
an algebraic skeleton. In the category of $LCAb$ of locally compact abelian topological groups
and continuous homomorphisms, with properly adjusted definitions, one can get a much more satisfying picture.

\smallskip

First af all, choose $\Cal{Z} := \bold{C}_1^*= \{ z\in \bold{C}^*\,|\, |z|=1\}$.
This is {\it a dualizing object}: for each $\Cal{K}$
in $LCAb$ there exists the internal $\roman{Hom}\, (\Cal{K},\Cal{Z})$
object, called the {\it character group} $\widehat{\Cal{K}}$,
and the map $\Cal{K} \mapsto \widehat{\Cal{K}}$ extends
to the equivalence of categories $LCAb \to LCAb^{op}$
(Pontryagin's duality).

\smallskip

Let now $\psi$ be a continuous cocycle $\Cal{K}\times\Cal{K}\to\Cal{Z}$
so that $\varepsilon$ is a continuous bicharacter.
Call the extension $\Cal{G} (\Cal{K}, \psi )$
a  {\it Heisenberg group}, if the map $x\mapsto \varepsilon (x,*)$ identifies $\Cal{K}$ with
$\widehat{\Cal{K}}$. In the case $\Cal{K}=\bold{R}^{2N}$, we call $\Cal{G}(\bold{R}^{2N},\psi)$ 
a {\it vector Heisenberg group}. 

\smallskip

Choose $k=\bold{C}$, and $\chi$ continuous. The formula (1.6) 
makes sense e.g. for continuous functions $f$.
Especially interesting, however, is the representation
on $L_2(\Cal{K})$ which makes sense  because the operators (21.6)
are unitary with respect to the squared norm $\int_{\Cal{K}} |\bold{f}|^2
d\mu_{Haar}$. Of course, square integrable functions
cannot be evaluated at points, so that $\bold{f}(x+y)$ in (1.6) must
be understood as the result of shifting $\bold{f}$ by $y\in \Cal{K}$;
similar precautions should be taken in the formula (1.8)
defining now the space $L_2(\Cal{K}//K_0)$ where $K_0$
is a closed isotropic subgroup (it is then automatically
liftable to a closed subgroup),
and in many intermediate calculations. See Mumford's treatment on
pp. 5--11 of [Mu] specially tailored for
readers with algebraic geometric sensibilities.

\smallskip

We will call $L_2(\Cal{K}//K_0)$, for a maximal isotropic subgroup $K_0\subset \Cal{K}$,
{\it  a basic representation} of the respective Heisenberg group.

\smallskip

The central fact of the representation theory of a
Heisenberg group $\Cal{G} (\Cal{K},\psi )$,
$\Cal{K}\in LCAb$, $\chi = id$, is this:

\medskip

{\bf 1.3.1. Theorem.} {\it (i) A basic representation of $\Cal{G} (\Cal{K},\psi )$ is unitary and irreducible.
 
\smallskip

(ii) Any unitary representation of 
 $\Cal{G} (\Cal{K},\psi )$ whose restriction
to the center is the multiplication by the identical character
is isomorphic to
the completed tensor product of $L_2(\Cal{K}//K_0)$
and a trivial representation. In particular, various basic represetations 
are isomorphic.}

\smallskip
If  (i), (ii) are satisfied, then the maps
$\Cal{Z}\to\Cal{G}:\,\lambda\mapsto (\lambda ,0),$
$\Cal{G}\to \Cal{K}:\, (\lambda ,x)\mapsto x$, 
describe $\Cal{G}$ as a central extension of
$\Cal{K}$ by $\Cal{Z}$:
$$
1\to \Cal{Z}\to \Cal{G}(\Cal{K},\psi )\to \Cal{K}\to 1.
$$
For a proof we refer the reader to [Mu].  
\bigskip
\centerline{\bf \S 2. Calculus of representations of Heisenberg groups.} 

\medskip

In this section we will recall some well--known facts about representations of Heisenberg groups. In the first part of this section we introduce 
the notions of matrix coefficients, integrability and square-integrablity of a representation on a Hilbert space. In the second part we will focus on the representations 
of vector Heisenberg groups.

\medskip

{\bf 2.1. Basic notions.} Let $\Cal{K}$ be an object of $LCAb$ and $\Cal{H}$ a Hilbert space. We 
consider a unitary strongly continuous irreducible representation $\rho$ of $\Cal{K}$ on $\Cal{H}$. For $\bold{f},\bold{g}\in\Cal{H}$ we call $S_\bold{g}\bold{f}(k)=\langle \bold{f},\rho(k)\bold{g}\rangle$ a {\it matrix coefficient} 
of the representation $\rho$ on $\Cal{H}$. If there exists a non-trivial $\bold{g}$ in $\Cal{K}$ such that $S_\bold{g}\bold{g}$ is square integrable with respect to the Haar measure on $\Cal{K}$, then the irreducible representation
$\rho$ is called {\it square integrable}. 

\smallskip

The irreducible representation $\rho$ is {\it integrable} if at least one matrix coefficient $S_\bold{g}\bold{g}$ is integrable for a non-trivial $\bold{g}\in\Cal{H}$, such a $\bold{g}$ is then  called an {\it integrable element}
of the irreducible and integrable representation $\rho$. We denote by $\Cal{A}_1(\Cal{K})$ the set of all integrable elements $g\in\Cal{H}$ and call it the {\it space of admissible atoms}.

\smallskip
  
The square integrability of $\rho$ implies the existence of {\it orthogonality relations} for matrix coefficients: for $\bold{f}_1,\bold{f}_2,\bold{g}_1,\bold{g}_2$ in $\Cal{H}$ the following holds, see e.g. [FeGr]:

$$
\int_{\Cal{K}}S_{\bold{g}_1}\bold{f}_1(k)\overline{S_{\bold{g}_2}\bold{f}_2(k)}dk=\langle\bold{f}_1,\bold{f}_2\rangle\overline{\langle \bold{g}_1,\bold{g}_2 \rangle} 
\eqno(2.1)
$$
Among important consequences of this identity we note the existence of a {\it reconstruction formula}
 for functions in $\Cal{H}$. Let $\bold{f}$ be in $\Cal{H}$ and $\langle\bold{g}_1,\bold{g}_2\rangle\ne 0$. Then we have 
$$
f= \langle\bold{g}_1,\bold{g}_2\rangle^{-1}\int_{\Cal{K}}S_{\bold{g}_1}\bold{f}_1(k)\rho(k)\bold{g}_2dk.
\eqno(2.2)
$$
In the case of the Heisenberg group $\Cal{G}(\Cal{K},\psi)$ this yields the most general version of the resolution of identity and the orthogonality relations are known as 
{\it Moyal's identity}. The space of admissible elements $\Cal{A}_1(\Cal{G(\Cal{K},\psi)})$ for the Heisenberg group $\Cal{G(\Cal{K},\psi)}$ provides a 
realization of {\it Feichtinger's algebra} \cite{Fe1}. Since the space of integrable elements of $\Cal{G(\Cal{K},\psi)}$ is another way of defining Feichtinger's algebra, see [FeGr]. Feichtinger's algebra is an important Banach space in harmonic anlaysis, 
time--frequency analysis and Gabor analysis. Later we discuss Feichtinger's algebra and some of its weighted versions in
the case of vector Heisenberg groups in more detail.
\smallskip
In our discussion of projective modules over quantum tori we also have to deal with smooth vectors of representations of a 
Heisenberg group $\Cal{G}(\Cal{K},\psi)$, where $\Cal{K}$ is an {\it elementary locally compact
abelian group}. Recall that this means $\Cal{K}$ is of the form {\it vector space $\times$ torus $\times$ 
lattice $\times$ finite group}.  We denote the Lie algebra of the Heisenberg group by $L$. A vector 
$\bold{f}\in \Cal{H}$ is called {\it smooth} if
for any $X_1, \dots ,X_n \in L$ the following expression
makes sense
$$
\delta U_{X_1}\circ \dots \delta U_{X_n} (\bold{f})
$$
where $\delta U_{X} (\bold{f})$ is defined as the limit when $t\to 0$
$$
\delta U_X(\bold{f}) := \roman{lim}\,\frac{U_{exp(tX)}\bold{f} - \bold{f}}{t}.
\eqno(2.3)
$$
It is known that the space $\Cal{H}_{\infty}$ of smooth vectors 
is dense and the operators $\delta U_X$ are skew adjoint but unbounded, [Mu].

\medskip

{\bf 2.2. Basic representations of vector Heisenberg groups.} Now we discuss the vector Heisenberg group $\Cal{G}(\bold{R}^{2N},\psi )$ in more detail. 
Note that in this case we can choose $\psi$ as an antisymmetric bicharacter with values in $\bold{C}^*_1.$ After choosing an appropriate
basis, we identify $\bold{R}^{2N}$ with the space of pairs of column vectors $x=(x_1,x_2),\, x_i\in \bold{R}^N$. In this case $\psi$ can be written in the form 
$$
\psi_A (x,y) = e^{\pi i A(x,y)}
\eqno(2.4)
$$
where $A:~ \Cal{K}\times \Cal{K} \to \bold{R}$
is a nondegenerate antisymmetric pairing,
such as
$$
A(x,y) = x_1^ty_2-x_2^ty_1
$$ 
where $x_i^t$ denotes the transposed row vector.
We have then $\varepsilon_A (x,y) =e^{2\pi iA(x,y)}$.

\smallskip

In particular, the subspace $\bold{R}^N\times\{0\}$ is a maximal
$\varepsilon$--isotropic closed subgroup of $\bold{R}^{2N}$.

\smallskip

We will recall the structure of two Heisenberg
representations of $\Cal{G}(\bold{R}^{2N},\psi_A )$ 
using normalizations adopted in [Mu].

\medskip

{\it Model I.} In this model the unitary representation of $\Cal{G}(\bold{R}^{2N},\psi )$ lives on $L_2(\bold{R}^N)$, the Hilbert space of square integrable
complex functions $\bold{f}$ on $\bold{R}^N$  with the scalar product
$$
\langle \bold{f},\bold{g}\rangle := \int \bold{f}(x_1)\overline{\bold{g}(x_1)}\,d\mu
\eqno(2.5)
$$ 
where $d\mu$ is the Haar measure in which $\bold{Z}^N$
has covolume 1.

\smallskip

The action of $\Cal{G} (\bold{R}^N,\psi_A )$, with central character
$\chi (\lambda ) =\lambda$, is given by the formula
$$
(U_{(\lambda ,y)}\bold{f}) (x_1)= \lambda\,
e^{2\pi i\,x_1^ty_2 +\pi i\,y_1^ty_2}\,\bold{f}(x_1+y_1).
\eqno(2.6)
$$
Many authors refer to the {\it Model I} as the {\it Schr\"{o}dinger representation} of $\Cal{G}(\bold{R}^{2N},\psi_A )$.
\smallskip
{\it Model II${}_T$.} The second model is actually
a family of models depending on the choice 
of a {\it compatible} K\"ahler structure
upon $\Cal{K}=\bold{R}^{2N}.$
A general K\"ahler structure on $\bold{R}^{2N}$
can be given by a pair consisting of a complex structure
and a Hermitean scalar product $H$. We will call this
K\"ahler structure compatible (with
the choice (2.4)) if $\roman{Im}\,H=A$. Such structures
are parametrized by the Siegel space $\Cal{S}_N$ consisting
of symmetric matrices $T\in M(N,\bold{C})$ with
positive definite $\roman{Im}\,T.$ 

\smallskip

In particular, the complex structure defined by $T$
identifies $\bold{R}^{2N}$ with $\bold{C}^N$
via 
$$
(x_1,x_2)=x \mapsto \underline{x} := Tx_1+x_2,
\eqno(2.7)
$$
and we have
$$
H_T(\underline{x},\underline{x})=\underline{x}^t\,(\roman{Im}\,T)^{-1}\,
\underline{x}^*
\eqno(2.8)
$$
where $*$ denotes the componentwise complex conjugation.

\smallskip

Consider the Hilbert space $\Cal{H}_T$ of holomorphic functions
on $\bold{C}^N=\Cal{K}$ consisting of the functions with finite
norm with respect to the scalar product
$$
\langle \bold{f}, \bold{g}\rangle_T := \int_{\bold{C}^N} \bold{f}(\underline{x})\,
\overline{\bold{g}(\underline{x})}\,
e^{-\pi H_T(\underline{x},\underline{x})}\,d\nu
\eqno(2.9)
$$
where $d\nu$ is the translation invariant measure making
$\bold{Z}^{2N}$ a lattice of covolume 1 in $\bold{R}^{2N}.$ 

\smallskip

For $(\lambda ,y)\in \Cal{G}(\bold{C}^N,\psi_A )$ and a holomorphic
function $f$ on $\bold{C}^N$, put
$$
(U^{\prime}_{(\lambda ,y)}\bold{f})(\underline{x}):=
\lambda^{-1}e^{-\pi H_T(\underline{x},\underline{y})
-\frac{\pi}{2} H_T(\underline{y},\underline{y}) }
\bold{f}(\underline{x}+\underline{y}).
\eqno (2.10)
$$
A straightforward check shows that these operators are 
unitary with respect to 
(2.9), and moreover, that they define a representation
of $\Cal{G}(\bold{C}^N,\psi_A )$ in $\Cal{H}_T$
corresponding to the character 
$\chi(\lambda )=\lambda^{-1}$ of $\bold{C}_1^*$,
in the sense of formula (1.8). This is generalization of the classical
Bargmann--Fock representation, which corresponds to the choice $T=\tau I$ for $\tau$ with a positive imaginary part. 

\smallskip

It turns out that this representation is irreducible on $\Cal{H}_T$
and thus is a model of the Heisenberg representation.

\smallskip

The proof of irreducibility spelled out in [Mu] 
involves constructing {\it vacuum vectors}
in $\Cal{H}_T$ which in this model turn out to be simply
constant functions. Translated via
canonical (antilinear) isomorphism into other models they look
differently. For example they become
(proportional to) a ``quadratic exponent''
$\bold{f}_T:=e^{\pi ix^t_1 Tx_1}$ in Model I ( i.~e. $L_2(\bold{R}^{2N}//\bold{R}^{N})$) 
or to an essentially classical theta--function
$e^{\pi ix_1^t \underline{x}}\vartheta (\underline{x},T)$
in $L_2(\bold{R}^{2N}//\bold{Z}^{2N})$.
They are called  ``theta--vectors'' in [Sch].
For details, see the Theorem 2.2 in [Mu] and the discussion
around it. 

\bigskip


\centerline{\bf \S 3. Principle notions of time--frequency analysis in the context }
\smallskip
\centerline{\bf  of vector Heisenberg groups.} 

\medskip
In this section we interpret  matrix coefficients of the Schr\"{o}dinger representation of the vector Heisenberg group in terms 
of time--frequency analysis. This leads us naturally to the study of time--frequency representations, namely the short--time Fourier transform and 
the Wigner distribution from the phase space approach to quantum mechanics.

\smallskip

{\bf 3.1. Time--frequency representations.} The basic task in signal analysis is to analyse the spectral content of a given signal $\bold{f}$, i.e. a complex--valued function $\bold{f}(t)$ on $\bold{R}^N$.  
Traditionally one invokes the {\it Fourier transform} $\Cal{F}\bold{f}$ to gain some insight about the 
spectral resolution of $\bold{f}$, i.e. the frequencies $\omega$ contained in the given signal $\bold{f}$. 
We define the Fourier transform $\Cal{F}\bold{f}$ for $\bold{f}\in L_1(\bold{R}^N)$ by 
$$
\Cal{F}\bold{f}(\omega)=\widehat{\bold{f}}(\omega)=\int_{\bold{R}^N}\bold{f}(t)e^{-2\pi i t\cdot\omega}dt.
\eqno(3.1)
$$
One drawback of this approach is that $\Cal{F}\bold{f}$ does not provide any local information about the signal $\bold{f}$. 
Therefore one is looking for joint time and frequency representations. In time--frequency analysis it is very useful to consider $\bold{R}^{2N}$ as 
$$
\bold{R^{2N}}=\bold{R}^N\times \widehat{\bold{R}}^N =\{(t,\omega )\,|\,t,\omega\in \bold{R}^N\}.
\eqno(3.2)
$$
If we want to emphasize this description of $\bold{R}^{2N}$ we refer to it as {\it time--frequency space}. 

\smallskip

In order to obtain information about the ``local frequency spectrum" of $\bold{f}$, we use a {\it window} $\bold{g}$ to localize the signal $\bold{f}$ and 
take the Fourier transform of this localized piece of $\bold{f}$. 

\smallskip

{\it Short--Time Fourier Transform (STFT)} of a signal $\bold{f}$ with respect to
{\it a window function $\bold{g}$}:
$$
V_\bold{g}\bold{f} (x,\omega ):=
\int_{\bold{R}^N} \bold{f}(t)\overline{\bold{g}} (t-x)\,e^{-2\pi i t\cdot\omega}dt.
\eqno(3.3)
$$
The properties of $V_\bold{g}\bold{f}$ depend crucially on the choice of the window $\bold{g}$. 
It turns out that the integrable and smooth vectors of the vector Heisenberg 
group provide good classes of window functions. If $\bold{f}$ and $\bold{g}$ are for example Schwartz functions on $\bold{R}^N$, 
then $V_\bold{g}\bold{f}(x,\omega)$ measures the amplitude 
of the frequency band near $\omega$ at time $x$. 

\smallskip

Among the many facts about STFT, we want to mention one on the relation between the 
STFT of $\bold{f}$ and the STFT of $\widehat{f}$. 
By an application of Parseval's theorem to $V_\bold{g}\bold{f}$ 
one obtains
$$
V_{\widehat{\bold{g}}}\widehat{\bold{f}}(x,\omega)=e^{-2\pi ix\cdot\omega}V_\bold{g}\bold{f}(\omega,-x).
\eqno(3.4)
$$
Therefore if we choose the Gaussian $g_0(t)=e^{-\pi t^2}$ as a window function (or any window invariant under Fourier transform), 
then the time--frequency content of $\widehat{\bold{f}}$ is just the one of $\bold{f}$ rotated
by an angle of $\pi/2$. 

\smallskip

During the last two decades a representation theoretic interpretation of the STFT has become of great 
relevance in time--frequency analysis, because it has allowed to put the time--frequency analysis on 
solid mathematical ground. The main proponents of this line of research are Feichtinger and Gr\"ochenig, 
who described in \cite{FeGr} a correct framework for a rigorous treatment of time--frequency analysis, 
the {\it coorbit theory}. Later we briefly discuss some aspects of their coorbit theory for the Schr\"{o}dinger 
representation of vector Heisenberg groups. 

\smallskip

In time--frequency analysis one associates to a point $(x,\omega)$ in the time--frequency plane a {\it time--frequency shift:} the unitary operator 
$\pi (x,\omega)$ on $L_2(\bold{R}^N)$ defined by
$$
\bold{f}\mapsto \pi (x,\omega )\bold{f}(t) :=e^{2\pi i t\cdot \omega} \bold{f}(t-x).
\eqno(3.5)
$$

The same operator occurs in the basic representation of the central extension
of the time--frequency space:
$$
\pi (x,\omega )= U_{(e^{-\pi it\cdot \omega},(x,\omega))}.
\eqno(3.6)
$$

The restriction of $\pi(x,\omega)$ to the maximal $\varepsilon_A$-isotropic subspaces $\bold{R}^N\times\{ 0\}$ and $\{0\}\times\bold{R}^N$ of the 
time--frequency space $\bold{R}^{2N}$ are respectively the {\it translation operator}
 $T_x:= \pi (x,0)$ and the {\it modulation operator}  $M_{\omega}:=\pi(0,\omega )$.

\smallskip 
The mapping $(x,\omega)\mapsto \pi(x,\omega)$ is a unitary (projective) strongly continuous representation of $\bold{R}^{2N}$ on $L_2(\bold{R}^N)$.
 Therefore, we can express $V_\bold{g}\bold{f}$ as the matrix-coefficient of this projective representation, i.e. 
 $$
 V_\bold{g}\bold{f}(x,\omega)=\langle\bold{f},\pi(x,\omega)\bold{g}\rangle.
 \eqno(3.7)
 $$
This intrinsic description of STFT in terms of the Heisenberg group amplifies the great relevance of the STFT and makes it a basic object of study. 
Moreover,  most other time--frequency representations, such as cross--Wigner distribution,  have a description in terms of the STFT.

\smallskip

Recall, that the {\it cross--Wigner distribution} of two signals $\bold{f}$ and $\bold{g}$ is defined by 
$$
  W(\bold{f},\bold{g})(x,\omega)=2^Ne^{4\pi ix\cdot\omega}\int_{\bold{R}^N} \bold{f}(x+\tfrac{t}{2})\overline{\bold{g}}(x-\tfrac{t}{2})e^{-2\pi it\cdot\omega} dt.
\eqno(3.8)
$$
The cross--Wigner distribution is just a short--time Fourier transform in disguise:
$$
   W(\bold{f},\bold{g})(x,\omega)=2^Ne^{4\pi ix\cot\omega}V_{\tilde{\bold{g}}}\bold{f}(2x,2\omega),
\eqno(3.9)
$$
where $\tilde{\bold{g}}(x)=\bold{g}(-x)$ denotes the reflection of $\bold{g}$ with respect to the origin.

\medskip

{\bf 3.2. Function spaces for time--frequency analysis.} Since STFT is one of the key players in the time--frequency analysis, 
it is natural to consider  function spaces defined  in terms of integrability conditions or decay conditions of the STFT. 
In the early 1980's, H.~G.~ Feichtinger introduced the so--called {\it modulation spaces} in exactly this manner [Fe2]. For a thorough 
discussion of the history and the relevance of modulation spaces in various branches of mathematics and engineering see the 
excellent survey article by Feichtinger [Fe3]. Among the class of modulation spaces one Banach space, {\it Feichtinger's algebra} $S_0(\bold{R}^N)$, 
stands out as the most prominent member, already introduced by different methods in [Fe1]. Feichtinger's algebra is a good 
substitute for the Schwartz class as long as one is not dealing with partial differential equations, see [FeKo]. Feichtinger's algebra $S_0(\bold{R}^N)$ 
has an intrinsic description in terms of the Heisenberg group, namely it is the space of integrable vectors for the Schr\"{o}dinger representation of 
the Heisenberg group $\Cal{G}(\bold{R}^{2N},\psi_A)$. 

\smallskip

Let $v_s$ be the submultiplicative weight $v_s(x,\omega)=(1+\|x\|^2_2+\|\omega\|^2_2)^{s/2}$ on $\bold{R}^{2N}$, 
i.e. $v_s(x+y,\omega+\eta)\le v_s(x,\omega)v_s(y,\eta)$ for $(x,\omega),(y,\eta)\in\bold{R}^{2N}$. Recall that a weight $m$ is called $v_s$-moderate, if 
$m(x+y,\omega+\eta)\le v_s(x,\omega)m(y,\eta)$ for all $(x,\omega),(y,\eta)\in\bold{R}^{2N}$. Let $g$ be a nonzero function in the Schwartz class $\Cal{S}(\bold{R}^N)$. 
Then a tempered distribution $\bold{f}\in\Cal{S}^\prime(\bold{R}^N)$ belongs to the {\it modulation space} $M_{p,q}^m(\bold{R}^N)$, if 
$$  
\|\bold{f}\|_{M_{p,q}^m}:=\Big(\int_{\bold{R}^N}\Big(\int_{\bold{R}^N}|V_\bold{g}\bold{f}(x,\omega)|^pm(x,\omega)^p\Big)^{q/p}d\omega\Big)^{1/q}<\infty.
\eqno(3.10)
$$
For $1\le p,q\le\infty$ the modulation spaces $M_{p,q}^m(\bold{R}^N)$ are Banach spaces, and different functions $\bold{g}$ yield equivalent norms on $M_{p,q}^m(\bold{R}^N)$.  We will write $M_p$ for $M_{p,p}.$
The properties of modulation spaces are by their definition related to the properties of the short--time Fourier transform. 

\smallskip

The modulation spaces 
$M_1^{v_s}(\bold{R}^N)$ are of great relevance in Gabor analysis, since they provide a natural class of window functions. We denote $M_1^{v_s}(\bold{R}^N)$ by $M_1^{s}(\bold{R}^N)$.
The space $M_1(\bold{R}^N)=M_1^{v_0}(\bold{R}^N)$ is the well-known Feichtinger algebra $S_0(\bold{R}^N)$ \cite{Fe3}.

\smallskip

Below we summarize some properties used later
 in our treatment of projective modules over 
quantum tori. 

\smallskip

{\bf 3.2.1. Proposition.} {\it  The following holds:

\smallskip

(i) $M_{p,q}^m(\bold{R}^N)$ is invariant under time--frequency shifts, i.e. 
$$
\|\pi(y,\eta)\bold{f}\|_{M_{p,q}^m}\le C v_s(y,\eta) \|\bold{f}\|_{M_{p,q}^m}.
\eqno(3.11)
$$

\smallskip

(ii) $M_{p,q}^m(\bold{R}^N)$ is invariant under Fourier transform, i.e. 
$$
\|\widehat{\bold{f}}\|_{M_{p,q}^m}\le C \|\bold{f}\|_{M_{p,q}^m}.
\eqno(3.12)
$$

\smallskip

(iii) Let $\bold{f},\bold{g}$ be in $M_1^{s}(\bold{R}^N)$. Then $V_\bold{g}\bold{f}$ is in $L_1(\bold{R}^N)$.}

\medskip

For proofs of these statements we refer to \cite{Gr2} or \cite{FeLu}. The basic strategy for proving  
(i) and (ii)  is to establish these properties for STFT. The statement (i) 
follows from the {\it covariance property} of STFT:

$$
V_{\bold{g}}(\pi(y,\eta)\bold{f})(x,\omega)=e^{2\pi iy\cdot\omega} V_{\bold{g}}\bold{f}(x-y,\omega-\eta). 
$$ 
The proof of (ii) relies on eq. (3.7). 

\smallskip

Now we want to give a global description of the smooth vectors $\Cal{H}_\infty$ for the vector Heisenberg group $\Cal{G}(\bold{R}^N,\psi_A)$. 
In this case $\Cal{H}_\infty$ is the space of Schwartz functions $\Cal{S}(\bold{R}^N)$. The basic observation is that 
$$
\Cal{S}(\bold{R}^N)=\bigcap_{s\ge 0}M_1^{s}(\bold{R}^N).
\eqno(3.13)
$$
Consequently, $f\in\Cal{S}(\bold{R}^N)$ if and only if $\|\bold{f}\|_{M_1^s}$ is finite for all $s\ge 0$, see \cite{Gr2} for a proof. Note 
that this description of $\Cal{S}(\bold{R}^N)$ allows one  to transfer many statements from $M_1^{s}(\bold{R}^N)$ to $\Cal{S}(\bold{R}^N)$. 

\smallskip
We want to close this section with a discussion of the {symplectic Fourier transform}. For $\bold{F}\in L_1(\bold{R}^{2N})$ we define 
its {\it symplectic Fourier transform} by 

$$
\Cal{F}_s\bold{F}(x,\omega)=\int_{\bold{R}^{2N}}\bold{F}(y,\eta)e^{2\pi i(y\cdot\omega-x\cdot\eta)}dyd\eta.
\eqno(3.14)
$$ 

Observe that for a fixed $(x,\omega)\in\bold{R}^{2N}$ the set $\{e^{2\pi i(y\cdot\omega-x\cdot\eta)}\,|\,(y,\eta)\in\bold{R}^{2N}\}$ 
is actually the character group of the time--frequency plane, i.e. $\widehat{\bold{R}}^N\times\bold{R}^N$. Therefore the 
symplectic Fourier transform is the ordinary Fourier transform of the time--frequency plane $\bold{R}^N\times\widehat{\bold{R}}^N$. 

\smallskip 

As in the Euclidan case,  $\bold{M_1}^{v_s\otimes v_s}(\bold{R}^{2N})$ is invariant under the symplectic Fourier transfrom and 
consequently $\Cal{S}(\bold{R}^{2N})$ is invariant under the symplectic Fourier transform. 
The following fact has important consequences in Gabor analysis  (Janssen representation of a Gabor frame operator,
see [FeKo], [FeLu], [Ja4]) and in 
the construction of projective modules over quantum tori (see [Ri5]).

\medskip

{\bf 3.2.2. Proposition.}  {\it Let $\bold{f}_1,\bold{f}_2,\bold{g}_1,\bold{g}_2$ be in $M_1^s(\bold{R}^N)$. Then }
$$
\Cal{F}_s(V_{\bold{g}_1}\bold{f}_1)\cdot \overline{(V_{\bold{g}_2}\bold{f}_2})(y,\eta)=\overline{(V_{\bold{f}_1}\bold{f}_2}\cdot V_{\bold{g}_2}\bold{g}_1)(y,\eta). 
\eqno(3.15)
$$ 

\medskip
In [FaSc] the authors point out that this identity has been known in the signal analysis community since the early 1960's, when Sussman obtained 
this fact in his work on time--varying filters \cite{Su}. Later Rieffel gave a proof for signals and windows in $\Cal{S}(\bold{R}^N)$ in [Ri5] , which follows from 
the statement as indicated above. In Gabor analysis Feichtinger and Kozek were the first who explicitly formulated this fact for Feichtinger's algebra $S_0(\bold{R}^N)$. 
In [FeLu] we have discussed this identity in great detail for signals in dual classes of modulation spaces. The main idea is to choose the signals in 
such a way that the product of the STFT's is in $M_1(\bold{R}^{2N})$ and then use its invariance under the symplectic Fourier transform to justify 
the application of the symplectic Fourier transform.


\bigskip
\bigskip
\centerline{\bf \S 4. Quantum tori associated to embedded lattices }
\smallskip
\centerline{\bf in the vector Heisenberg groups.}

\medskip

We briefly recall the notions of quantum theta functions, Heisenberg group of quantum tori and quantum tori 
associated to embedded lattices. Our presentation follows largely the lines of \cite{Ma5}. 

\medskip
{\bf 4.1. Heisenberg group of quantum torus and quantum theta functions.}
Let $H$ be a free abelian
group of finite rank written additively, $k$ a ground field, and $\alpha :\,H\times H\to 
k^*$ a skewsymmetric pairing.
The {\it quantum torus} $T(H,\alpha )$
with the character group $H$ and quantization
parameter $\alpha$ is represented by an algebra
generated by a family of formal exponents $\{e(h)=e_{H,\alpha}(h),\,h\in H\}$,
satisfying the relations
$$
e(g) e(h)= \alpha (g,h) e(g+h).
\eqno(4.1)
$$
In particular, $T(H,1)$ is an algebraic torus, spectrum
of the group algebra $k[e(h)\,|\, h\in H]$ of $H$. The group of its
points $x\in T(H,1)(k) = \roman{Hom}\, (H,k^*)$ acts upon
functions on $T(H,\alpha )$ mapping $e_{H,\alpha}(h)$
to $x^*(e_{H,\alpha}(h)):=h(x) e_{H,\alpha}(h)$ where $h(x)$
denotes the value of the character $e(h)$ at $x$. 
 
\smallskip

The {\it Heisenberg group}  of the quantum torus $T(H,\alpha )$ introduced
in [Ma3] and denoted there $\Cal{G} (H,\alpha )$ consists of   
all maps of the form
$$
\Phi \mapsto c\, e_{H,\alpha}(g)\, x^*(\Phi )\, e_{H,\alpha}(h)^{-1},\
c\in k^*;\, x\in T(H,1)(k);\, g,h\in H,
\eqno(4.2)
$$
where $\Phi$ is a formal theta function, see [Ma3] Eq.(0.19). Any such map has a 
unique representative of the same form
in which $h=0$ (``left representative''). Writing this
representative as $[c; x, g]$ we get the composition law
$$
[c^{\prime};x^{\prime},g^{\prime}][c;x,g]=
[c^{\prime}c\, g(x^{\prime})\, \alpha(g^{\prime},g);x^{\prime}x,g^{\prime}+g].
\eqno(4.3)
$$
In other words, this group is the central extension of
$\roman{Hom}\,(H,k^*)\times H$ by $k^*$ corresponding
to the bicharacter
$$
\psi ((x^{\prime},g^{\prime}),(x,g))=g(x^{\prime}) \alpha (g^{\prime},g)
\eqno(4.4)
$$
and having the associated bicharacter
$$
\varepsilon ((x^{\prime},g^{\prime}),(x,g))=g(x^{\prime})
g^{\prime}(x)^{-1} \alpha^2(g^{\prime},g).
\eqno(4.5)
$$
In particular, if a subgroup $B\subset \roman{Hom}\,(H,k^*)\times H$
is liftable to $\Cal{G} (H,\alpha )$, the form (4.4)
restricted to $B$ must be symmetric: this is the main part
of Lemma 2.2 in [Ma3]. 

\smallskip

A lift $\Cal{L}$ of $B$ to a subgroup of $\Cal{G} (H,\alpha )$ is called
{\it a multiplier}. The  restriction to $B$ of the form (4.4),
$\langle \,,\rangle :\, B\times B\to k^*$, is called {\it the structure
form} of this multiplier. (Formal) linear combinations of the
exponents $e_{H,\alpha}$ invariant with respect to the
action of $\Cal{L} (B)$ constitute a linear space $\Gamma (\Cal{L})$
and are called (formal) {\it quantum theta functions}. 

\medskip

{\bf 4.2. Embedded lattices and quantum tori.} In this section $\Cal{K}$ denotes 
an object of $LCAb$, $\psi$ is a
bicharacter of $\Cal{K}$ such that 
$\varepsilon$ (cf. (1.5)) identifies $\Cal{K}$ with
$\widehat{\Cal{K}}$. Let $\Cal{G} (\Cal{K},\psi )$ be
the respective Heisenberg group, central extension of
$\Cal{K}$ by $\Cal{Z} = \bold{C}_1^*$ as above.

\smallskip

We will call {\it an embedded lattice} 
a closed subgroup $D\subset \Cal{K}$ such that
$D$ is a finitely generated free abelian group,
whereas $\Cal{K}/D$ is a topological torus. In this section we
consider only those groups $\Cal{K}$ which admit embedded lattices. 

\smallskip

Consider a family of constants 
$$
\{c_h\in \bold{C}_1^*,h\in D\}.
$$ 
Put 
$$
E(h):= (c_h,h)\in \Cal{G} (\Cal{K},\psi )
\eqno(4.6)
$$ 
From (1.2) we get
$$
E(g)E(h)=\frac{c_gc_h}{c_{g+h}}\,\psi (g,h)\,E(g+h).
$$
Assume that
$$
\alpha (g,h) := \frac{c_gc_h}{c_{g+h}}\,\psi (g,h)
\eqno(4.7)
$$
is a skewsymmetric pairing. Then the map $e_{D,\alpha}(h)
\mapsto E(h)$ is compatible with the relations (1.6), i.e. 
defines a cohomologous representation of the Heisenberg group. 
In particular any representation $U$ of
$\Cal{G} (\Cal{K},\psi )$ induces a representation of
an appropriate function algebra of the quantum torus
$T(H,\alpha )$. One easily sees that any $\alpha$ on $D$ can be induced
from an appropriate lattice embedding of $D$;
one can even take $\psi$ to be a skewsymmetric bicharacter so that
$\alpha$ will coincide with the restriction of $\psi$.

\smallskip

The condition (4.1) in the definition of formal exponents may be 
considered as a {\it projective representation} of $D$ with respect to 
the 2-cocycle $\alpha$. There is a canoncial correspondence between these 
representations and involutive representations of the {\it twisted group algebra} 
$\ell_1(D,\alpha)$ of $D$. It is known as the method of integrated representations.

\smallskip

The {\it twisted group algebra} $\ell_1(D,\alpha)$ of $D$ consists of all absolutely 
summable sequences $\bold{a}=(a_h)_{h\in D}$ where the multiplication is defined 
as {\it twisted convolution} of $\bold{a}$ and $\bold{b}$ by 
$$ 
\bold{a}\natural\bold{b}(h)=\sum_{l\in D}a_{l}b_{h-l}\alpha(h-l,l),
\eqno(4.8)
$$
and involution $\bold{a}^*=(a^*_h)_{h\in D}$ of $\bold{a}$ is defined by  
$$
a_h^*=\overline{\alpha(h,h)a_{-h}}.
\eqno(4.9)
$$ 
Consequentely $\ell_1(D,\alpha)$ becomes a Banach algebra with respect to the norm 
$\|\bold{a}\|_1$.

\smallskip

If we ``integrate" the formal exponents $\{e_h:h\in H\}$, then we get an {\it involutive representation} of 
$\ell_1(D,\alpha)$ as follows: 
$$
\bold{a}=(a_h)_{h\in D}\mapsto \rho_D(\bold{a}):=\sum_{h\in D}a_he_{D,\alpha}(h).
\eqno(4.10)
$$
More explicitely, we have for all $\bold{a},\bold{b}$ in $\ell_1(D)$: 
$$
\rho_D(\bold{a})\rho_D(\bold{b})=\rho_D(\bold{a}\natural\bold{b})~~\text{and}~~\rho_D(\bold{a})^*=\rho_D(\bold{a}^*).
\eqno(4.11)
$$
It is a well--known fact that there is a one--to--one correspondence between projective representations of $D$ and 
involutive representations of the twisted group algebra $\ell_1(D,\alpha)$. Finally we consider the {\it twisted group} 
$C^*$--algebra $C^{*}(D,\alpha )$ of $D$, which is the enveloping $C^*$-algebra of $\ell_1(D,\alpha)$.

\smallskip

Later we want to construct projective modules over {\it smooth subalgebras} $\Cal{A}$ of $C^{*}(D,\alpha )$ in the sense of A.~Connes,. This means that
 $\Cal{A}$ is stable under the holomorphic function calculus of $C^{*}(D,\alpha )$. The algebra $C^{\infty}(D,\alpha )$ of 
elements $\sum_{h\in D}a_he_{D,\alpha}(h)$ with coefficients $\{ a_h\}_{h\in D}$ belonging to the Schwartz space $\Cal{S}(D)$ is 
a well-known example of a smooth subalgebra of $C^{*}(D,\alpha )$. 

\smallskip

We want to point out that investigations in signal analysis have given rise to 
a whole class of smooth subalgebras of $C^{*}(D,\alpha )$, cf.  [GrLe]. 
They are denoted $C_1^s(D,\alpha)$, 
where one imposes on the elements $\sum_{h\in D}a_he_{D,\alpha}(h)$ of $C^{*}(D,\alpha )$ the following summability conditions:
$$
\sum_{h\in D}|a_h|(1+|h|^2)^{s/2}<\infty.
\eqno(4.12)
$$
The fact that $C_1^s(D,\alpha)$ are smooth subalgebras of $C^{*}(D,\alpha )$ was shown by Gr\"ochenig and Leinert in [GrLe], where 
they proved the so--called {irrationality conjecture} of Gabor analysis. Recently J.~Rosenberg has given another proof 
for the case $s=0$ of the theorem of Gr\"ochenig and Leinert in [Ro].

\smallskip

Notice that $C^{\infty}(D,\alpha )=\bigcap_{s\ge 0}C_1^s(D,\alpha)$. In other words,
 one might consider $C_1^s(D,\alpha)$ as noncommutative analogues 
of differentable functions of order $s$. In [Lu2], we have explored this point of view in detail. 

\smallskip 

Alternatively, any element of $C^{*}(D,\alpha )$ can also be written
as a formal series  $\sum_{h\in D}a_he_{D,\alpha}(h)$
but there is no transparent way to specify
which sequences $\{ a_h\in \bold{C}\,|\,h\in D\}$
can occur as their ``noncommutative Fourier coefficients''. In this picture $C_1(D,\alpha)$ are noncommutative analogs of Wiener's algebra 
of Fourier series with absolutely convergent Fourier series. Recently Rosenberg wrote a very interesting paper \cite{Ro} stressing this point of view. 
The theory of Gabor frames sheds some light on the structure of the noncommutative Fourier coefficients.
It shows in particular  that the class of modulation spaces 
has a natural characterization in terms of such noncommutative Fourier coeffiecients. We come back to this issue in a later section.

\medskip

{\bf 4.3. Projective modules over smooth subalgebras of quantum tori.} In this section we discuss Rieffel's projective modules over $C^\infty(D,\alpha)$ and 
their extension to the setting of noncommutative Wiener algebras $C_1^s(D,\alpha)$, where $D$ is an embedded lattice of $\Cal{K}$. Our presentation 
is inspired by the results in [Lu1], [Lu3], [Lu4], [Ma5], [Ri5].  

\smallskip
In [Ri5], it is shown that the space of smooth vectors $\Cal{H}_\infty$ gives rise to a finitely generated projective $C^{\infty}(D,\alpha)$-module. Now 
we want to formulate the results of [Lu3], [Lu4] in this general framework. Therefore we introduce a family of subspaces of {\it admissible elements} $\Cal{A}_1^s(\Cal{K})$ for
the Heisenberg group $\Cal{G}(\Cal{K},\psi)$, as those elements $\bold{g}$ of $L_2(\Cal{K})$ with $V_{\bold{g}}\bold{g}$ in $L_1^s(\Cal{K})$. 
Note that $\Cal{H}_\infty=\bigcap_{s\ge 0}\Cal{A}_1^s(\Cal{K})$. 

\smallskip

For $\Phi,\Psi$ in $\Cal{A}_1^s(\Cal{K})$ 

$$
{}_D\langle \Phi ,\Psi\rangle :=
\sum_{h\in D} \langle \Phi , e_{D, \alpha}(h) \Psi\rangle\,e_{D,\alpha}(h).
\eqno(4.13)
$$
defines a scalar product with values in $C_1^s(D,\alpha)$, which is bounded on $\Cal{A}_1^s(\Cal{K})$. The space $\Cal{A}_1^s(\Cal{K})$ becomes 
a left pre--inner product $C_1^s(D,\alpha)$--module with respect to the following left action:  

$$
\bold{a}\cdot\bold{f}=\sum_{h\in D} a_he_{D,\alpha}(h)\cdot\bold{f},
\eqno(4.14)
$$
for $\bold{a}\in\ell_1^s(D)$ and $\bold{f}\in\Cal{A}_1^s(\Cal{K})$.

\smallskip

More explicitely, this means the following:

\smallskip

(i) {\it Symmetry:} ${}_D\langle \Phi ,\Psi\rangle^*=
{}_D\langle \Psi ,\Phi\rangle .$

\smallskip

(ii) {\it (Bi)linearity:} ${}_D\langle a\Phi ,\Psi\rangle =
a\,{}_D\langle \Phi ,\Psi\rangle$ for any $a\in C^{\infty}(D,\alpha )$.

\smallskip

(iii) {\it Positivity:}  ${}_D\langle \Phi ,\Phi\rangle$ 
belongs to the cone of positive elements of $C^*(D,\alpha)$.
Moreover, if ${}_D\langle \Phi ,\Phi\rangle =0$ then $\Phi =0.$

\smallskip

(iv) {\it Density:} The image of ${}_D\langle\,,\rangle$ is dense
in $C_1^s(D,\alpha )$.  

\medskip
Consequently the completion of $\Cal{A}_1^s(\Cal{K})$ with respect to $\|\Phi\|:=\|{}_D\langle \Phi ,\Phi\rangle\|^{1/2}$ becomes 
a finitely generated left $C^*(D,\alpha)$-module $P$. Since $C_1^s(D,\alpha)$ is a smooth subalgebra of $C^*(D,\alpha)$ the finitely generated 
left $C_1^s(D,\alpha)$--module $P_1^s$ is isomorphic to $P$, see Lemma 4 on p.~52 of  [Co2] and the discussion around Proposition 3.7 in [Ri5]. 

\smallskip

As indicated above, this implies that $\Cal{H}_\infty$ is also a projective finitely generated left $C^\infty(D,\alpha)$-module.
\medskip 

{\bf 4.4. Projective modules over quantum tori for dual embedded lattices.} Let $D\subset \Cal{K}$ be an embedded
lattice as in 4.2. Denote by $D^!$ the maximal closed subgroup
of $\Cal{K}$ orthogonal to $D$ with respect to $\varepsilon$.
From the Pontryagin duality it follows that $D^!$ (resp. $D$) can be
canonically identified with the character group of 
$\Cal{K}/D$ (resp. $\Cal{K}/D^!$) so that $D^!$ is an embedded lattice as well.

\smallskip

Assume moreover that $\psi$ is an antisymmetric pairing, so that one
can choose $E(h) = (1,h)\in \Cal{G} (\Cal{K}, \psi )$ for $h\in D$
and for $h\in D^!$ and 
define on $\Cal{A}_1^s(\Cal{K})$ the structure of $C_1^s(D^!,\alpha^!)$--module
as well where $\alpha^!$ is the pairing induced
on $D^!$ by $\psi$. Operators $e_{D,\alpha}(h),\,h\in D,$ commute 
with operators $e_{D^!,\alpha^!}(g),\,g\in D^!.$ 

\smallskip

We can consider $\Cal{A}_1^s(\Cal{K})$ as a right 
$C_1^s (D^!,\alpha^!)^{op}$--module. Moreover, we can and will 
identify the latter algebra with $C_1^s (D^!,\overline{\alpha}^!)$
by $e_{D^!,\alpha^!}(h) \mapsto e_{D^!,\overline{\alpha}^!}(h)^{-1}$
and extending this map by linearity.

\medskip 

{\bf 4.5. Theorem.} {\it (i) We have 
$\|{}_D\langle \Phi ,\Phi\rangle\|^{1/2}
=\|{}_{D^!}\langle \Phi ,\Phi\rangle\|^{1/2}$.
The completion $\Cal{H}$ of $\Cal{A}_1^s(\Cal{K})$
with respect to this norm is a projective left
module over both quantum tori $C^{*}(D,\alpha )$ and $C^{*}(D^!,{\alpha}^!)$,
and each of these algebras is a total commutator of the other one.

\smallskip

(ii) Let  $C^{*}(D^!,\overline{\alpha}^!)$ 
act upon $\Cal{H}$ on the right as explained above. Consider the analog
of the scalar product (4.13)
$$
\langle \Phi ,\Psi\rangle_{D^!} :=\frac{1}{\roman{vol}\,\Cal{K}/D}\,
\sum_{h\in D} \langle  e_{D^!, \alpha^!}(h) \Psi ,\Phi 
\rangle\,e_{D^!,\overline{\alpha}^!}(h) \, \in C^{*}(D^!,\overline{\alpha}^!)
\eqno(4.15)
$$
It satisfies relations similar to (i)--(iv), and moreover,
for any $\Phi ,\Psi , \Xi$ the following
associativity relation holds:
$$
{}_D\langle \Phi ,\Psi\rangle \Xi = \Phi\,\langle \Psi ,\Xi\rangle_{D^!} .
\eqno(4.16)
$$
}
\smallskip
Consequentely this result holds also for the smooth vectors $\Cal{H}_\infty$ and the smooth subalgebra $C^\infty(D,\alpha)$. We discuss 
this theorem in the framework of Gabor analysis in a later section. 

\bigskip
\centerline{\bf \S 5. Quantum tori for embedded lattices and Gabor analysis.}
\medskip

In the first part of this section we introduce the basic notions of Gabor analysis and in the second part 
we use this framework to interpret projective modules over quantum tori in terms of Gabor frames.

\medskip

{\bf 5.1. Basics of Gabor analysis.} Recall that the square--integrability of the Schr\"odinger representation of the time--frequency plane yields 
the existence of a reconstruction formula for each $f\in L_2(\bold{R}^N)$:
$$
\bold{f}=\langle\bold{h},\bold{g}\rangle^{-1}\iint_{\bold{R}^{2N}}V_{\bold{g}}\bold{f}(x,\omega)\pi(x,\omega)\bold{h}~dxd\omega,
\eqno(5.1)
$$
for $\bold{g},\bold{h}\in L_2(\bold{R}^N)$ with $\langle\bold{g},\bold{h}\rangle\ne 0$. 

\smallskip

In the reconstruction formula (5.1) the time--frequency content of a signal $\bold{f}$ is analysed with respect to the system 
$\{\pi(x,\omega)\bold{g}:(x,\omega)\in\bold{R}^N\}$ for a window $\bold{g}\in L_2(\bold{R}^N)$, i.e. one considers the  
STFT $(\langle\bold{f},\pi(x,\omega)\bold{g}\rangle:(x,\omega)\in\bold{R}^N)$ for each building block $\pi(x,\omega)\bold{g}$, and then this information about the signal 
is used in the synthesis process with respect to the system $\{\pi(x,\omega)\bold{h}:(x,\omega)\in\bold{R}^N\}$. We call the function $\bold{h}$ the {\it synthesis window}. 
Note that there is just one restriction on the synthesis window, namely  $\bold{h}$ has to be non--orthogonal to the window $\bold{g}$. The reconstruction formula (5.1) 
is  well--known  in quantum mechanics, where it is refered to as {\it resolution of identity}. 

\smallskip  

In $1946,$ D. Gabor was looking for an ``optimal" way to transmit a signal ([Ga]). Therefore he wanted to maximize the content of information gained from the analysis process and 
to minimize the synthesis process. First, he suggested to use well--localized  window functions. Since the Gaussian $\bold{g}_0(t)=e^{-\pi t^2}$ 
(and its time--frequency shifts) is the minimizer of the Heisenberg uncertainty principle, Gabor relied his investigations on the Gaussian as window function. Second,  
Gabor  considered discrete analogues of the resolution of identity, where he replaced the time--frequency plane $\bold{R}^{2N}$ by the lattice $\bold{Z}^{2N}$. Relying upon an 
heuristic argument Gabor claimed that each $\bold{f}\in L_2(\bold{R}^N)$ has an expansion of the following type
$$
  \bold{f}=\sum_{k,m\in\bold{Z}^{N}}\langle \bold{f},\pi(k,m)\bold{g}_0\rangle\pi(k,m)\bold{h},
\eqno(5.2)  
$$
for some $\bold{h}$ in $L_2(\bold{R}^N)$. 

\smallskip

Until the late 1970s Gabor's paper [Ga] drew little attention of engineers and mathematicians, because the actual implementation of Gabor's expansions did 
not perform very well. In a series of papers [Ja1], [Ja2], [Ja3] a  mathematician Janssen undertook a rigorous investigation of Gabor's original expansions (5.2). 
The main result in [Ja1] shows that the series (5.2) converges for $\bold{f}$ and $\bold{g}$ not in $L_2(\bold{R}^N)$, because the coefficients turn out to grow logarithmically. 
Janssen proved instead convergence in the sense of tempered distributions and thereby 
explained why the expansions (5.2) are numerically unstable. The main reason for the problems with Gabor's original proposal is that the corresponding system does not give a frame for $L^2(\bold{R}^N)$. 

\smallskip

After the rigorous analysis of (5.2) mathematicians and engineers started looking for systems of functions in $L_2(\bold{R}^N)$ that would allow one to 
get numerically stable expansions of Gabor's type. The great breakthrough in this direction was the work [DaGrMe] in 1986, where the authors demonstrated that 
so--called {\it frames} of a Hilbert space provide stable reconstruction formulas. 

\smallskip

The notion of frames of a Hilbert space had been already introduced by Duffin and Schaeffer in their work on non--harmonic Fourier series [DuSc]. Let $\Cal{H}$ be 
a separable Hilbert space. Then a system $\{\bold{g}_i\,|\,i\in I\}$ is a {\it frame} for $\Cal{H}$ if and only if for each $\bold{f}\in\Cal{H}$ the following holds: there exist 
finite positive constants $A,B$ such that 
$$
  A\|\bold{f}\|^2_{\Cal{H}}\le\sum_{i\in I}|\langle \bold{f},\bold{g}_i\rangle|^2\le B\|\bold{f}\|^2_{\Cal{H}}.
\eqno(5.3)
$$ 
In \cite{DaGrMe} the authors investigated two kinds of frames generated by (i) Gabor systems and (ii) wavelets. Both systems are so--called {\it atomic decomposition}, 
because the elements of the systems are generated out of a building block $\bold{g}$ by the action of a group representation. A relevant group is the Heisenberg group in 
the case of Gabor frames and the affine group for wavelet systems. Shortly after the groundbreaking work of Daubechies, Grossman and Meyer, the mathematicians 
Feichtinger and Gr\"ochenig realized that Gabor systems and wavelets are just two examples of a general framework, which culminated in their {\it coorbit theory} [FeGr] 
and revealed the great relevance of representation theory for the construction of frames. 

\smallskip

After this brief historical review we continue our discussion of Gabor frames. Since the work of Feichtinger and Gr\"ochenig, a {\it Gabor system} 
$\bold{G}(\bold{g},D)$ is defined as the set of time--frequency shifts of a given {\it Gabor atom} $\bold{g}$ for an embedded lattice $D$ in $\bold{R}^{2N}$, i.e. 
$$
\bold{G}(\bold{g},D)=\{\pi(h)\bold{g}:h\in D\}.
\eqno(5.4)
$$  
A Gabor system $\bold{G}(\bold{g},D)$ is called a {\it Gabor frame} for $L_2(\bold{R}^N)$ if 
satisfies (5.3) for some constants $A,B$. 

\smallskip

The main task of Gabor analysis is to find reconstruction formulas for a function $\bold{f}$ in terms of $\bold{G}(\bold{g},D)$. 

\smallskip

The following operators associated to a Gabor system $\bold{G}(g,\Lambda)$  allow one to write down such reconstruction formulas. 

\smallskip

(i) The {\it analysis operator} $C_{\bold{g},D}$ maps functions $\bold{f}\in L_2(\bold{R}^N)$ to sequences on $D$ by 
$$
  C_{\bold{g},D}\bold{f}(h)=\langle \bold{f},\pi(h)\bold{g}\rangle\,h\in D.
\eqno(5.5)  	
$$

(ii) The {\it synthesis operator} maps sequences ${\bold{c}}=(c_h)$ on $D$ to functions on $\bold{R}^N$ as follows:
$$
{\bold{c}}\mapsto\sum_{h\in D}c_h\pi(h)\bold{g}.
\eqno(5.6)
$$

The {\it Gabor frame operator} $S_{\bold{g},D}$ corresponding to the Gabor system $\bold{G}(\bold{g},D)$ maps a function $\bold{f}$ to 
$$
  S_{\bold{g},D}\bold{f}=\sum_{h\in D}\langle \bold{f},\pi(h)\bold{g}\rangle\pi(h)\bold{g}.
\eqno(5.7)  
$$

Observe that $S_{\bold{g},D}$ is the composition of the analysis operator followed by the synthesis operator. Observe that $\bold{G}(\bold{g},D)$ is a Gabor frame for 
$L_2(\bold{R}^N)$ if and only if the Gabor frame operator $S_{\bold{g},D}$ is $\underline{invertible}$ on $L_2(\bold{R}^N)$. 

\smallskip

Now, the existence of reconstruction formulas for $\bold{f}$ is linked to the invertibility of the Gabor frame operator $S_{\bold{g},D}$. 
We define the {\it canonical dual Gabor atom} and the {\it canonical tight Gabor atom} of the Gabor frame $\bold{G}(\bold{g},D)$ by 
$\tilde{\bold{g}}:=S^{-1}_{\bold{g},D}\bold{g}$ and $\tilde{\bold{h}}=S^{-1/2}_{\bold{g},D}\bold{g}$, respectively. 
Then we have the following reconstruction formulas for $\bold{f}\in L_2(\bold{R}^N)$:
$$
  \bold{f}=\sum_{h\in D}\langle \bold{f},\pi(h)\bold{g}\rangle\pi(h)\tilde{\bold{g}}
\eqno(5.8)
$$
and
$$
 \bold{f}=\sum_{h\in D}\langle \bold{f},\pi(h)\tilde{\bold{h}}\rangle\pi(h)\tilde{\bold{h}}.
\eqno(5.9)
$$
For $\bold{g},\bold{h}\in L_2(\bold{R}^N)$ we call 
$$
  S_{\bold{g},\bold{h},D}\bold{f}=\sum_{h\in D}\langle \bold{f},\pi(h)\bold{h}\rangle\pi(h)\bold{g}
\eqno(5.11)
$$
a {\it Gabor frame-type operator}. These operators appear naturally in reconstruction formulas, 
e.g if $\tilde{\bold{g}}=S^{-1}_{\bold{g},D}\bold{g}$ then $S_{\bold{g},\tilde{\bold{g}},D}={I}$.

\smallskip

Around 1995 the papers [DaLaLa], [Ja4], [RoSh] developed the so-called {\it duality theory of Gabor analysis}, 
which marked a turning point in Gabor analysis. At the heart of all these contributions is the crucial observation, 
that a Gabor frame operator has a decomposition with respect to another Gabor system, the so-called 
{\it Janssen representation} of the Gabor frame operator. This relies on the observation that the 
Gabor frame operator $S_{\bold{g},D}$ commutes with $\pi(h)$ for $h\in D$:
$$
  \pi(h)S_{\bold{g},D}=S_{\bold{g},D}\pi(h)~~h\in D.
\eqno(5.12)  
$$
Therefore, it is natural to consider the lattice $D^!$ consisting of all points of $\bold{R}^{2N}$ commuting with $\{\pi(h):h\in D\}$:
$$
  D^!=\{h^!\in\bold{R}^{2N}:\pi(h)\pi(h^!)=\pi(h^!)\pi(h)~~\text{for all}~~h\in D\}.
\eqno(5.13)
$$

The lattice $D^!$ is the so-called {\it adjoint lattice}, which had been introduced by 
Feichtinger and Kozek ([FeKo]) in their discussion of the Janssen representation. In Gabor analysis the 
adjoint lattice is usually denoted by $D^\circ$. In other contexts this object has appeared in the work of Mumford on theta functions and polarized abelian varieties, Rieffel's construction of projective modules over noncommutative tori \cite{Ri5}. Note that the dual embedded lattice in the previous section is the adjoint lattice for the Heisenberg group $\bold{G}(\bold{R}^{2N},\psi_{I})$.

\smallskip

Consider $\bold{g},\bold{h}\in L_2(\bold{R}^{N})$ satisfying the following  condition of Tolimieri and Orr:
$$
  \sum_{h^!\in D^!}|\langle \bold{g},\pi(h^!)\bold{h}\rangle |<\infty.
\eqno(5.14)
$$
Then the Gabor frame-type operator $S_{g,h,D}$ has the following {\it Janssen representation}:
$$
  S_{\bold{g},\bold{h},D}f={\text{vol}(D)}^{-1}\sum_{h^!\in D^!}\langle \bold{h},\pi(h^!)\bold{g}\rangle\pi(h^!)\bold{f},
\eqno(5.15)
$$
where $\text{vol}(D)$ denotes the volume of a fundamental domain of $D$.

\smallskip

The Janssen representation follows from the {\it Fundamental Identity of Gabor analysis} (FIGA) after an application of a {\it symplectic Poisson summation formula}, see Eq. (5.18) below. The validity 
of Poisson summation formulas is a delicate matter. In the present situation the space of admissible vectors of the vector Heisenberg 
group $\Cal{G}(\bold{R}^{2N},\psi_A)$ provides  a good class of functions. Due to its great importance in harmonic analysis, time--frequency 
analysis and approximation theory this space \cite{Fe1} is known as {\it Feichtinger's algebra} $S_0(\bold{R}^N)$. We introduce weighted variants of Feichtinger's algebra $M_1^s(\bold{R}^N)$, 
these are elements of the class of {\it modulation spaces} \cite{Fe1}. The tempered distribution 
$\bold{f}\in\Cal{S}^{\prime}(\bold{R}^N)$ is in $M_1^s(\bold{R}^N)$ if 

$$
\|\bold{f}\|_{M_1^s}:=\iint_{\bold{R}^{2N}}|V_{\bold{g}}\bold{f}(x,\omega)|(1+|x^2|)^{s/2}dxd\omega<\infty
\eqno(5.16)
$$
for a window $\bold{g}$ in $\Cal{S}(\bold{R}^N)$. The defintion of $M_1^s(\bold{R}^N)$ is independent of the window function 
$\bold{g}$. It is a Banach space invariant under time--frequency shifts. The space $M_1(\bold{R}^N)$ is Feichtinger's algebra $S_0(\bold{R}^N)$, which 
is the minimal element in the class of Banach spaces invariant under time-frequency shifts. 

\smallskip

We formulate the symplectic Poisson summation formula in the following statement:

\smallskip

{\bf 5.1.1. Propostion.} {\it Let $D$ be an embedded lattice and $D^!$ its adjoint lattice in $\bold{R}^{2N}$. Then for $\bold{F}$ in $M_1(\bold{R}^{2N})$  
or $\bold{F}$ in $\Cal{S}(\bold{R}^N)$ we have 
$$
\sum_{h\in D}\bold{F}(h)=\frac{1}{vol(D)}\sum_{h^!\in D^!}\bold{F}(h^!).
\eqno(5.17)
$$
This holds pointwise, and the convegence is absolute.}

\smallskip

Therefore we are able to obtain the FIGA, which is a consequence of the 
symplectic Poisson summation formula for $\bold{F}=V_{\bold{g}_1}\bold{f}_1\cdot V_{\bold{g}_2}\bold{f}_2$. We just have to find conditions 
on the windows and functions guaranteeing that the product of the two STFT's is in $M_1(\bold{R}^N)$ or in $\Cal{S}(\bold{R}^N)$.
\medskip

{\bf 5.1.2. Propostion.} {\it Let $\bold{f}_1,\bold{f}_2,\bold{g}_1,\bold{g}_2$ be in $M_1(\bold{R}^N)$ or in $\Cal{S}(\bold{R}^N)$. Then we have }
$$
\sum_{h\in D}(V_{\bold{g}_1}\bold{f}_1 V_{\bold{g}_2}\bold{f}_2)(h)=\frac{1}{vol(D)}\sum_{h^!\in D^!}(V_{\bold{g}_1}\bold{g}_2 V_{\bold{f}_1}\bold{f}_2)(h^!). 
\eqno(5.18)
$$

The Janssen representation of $S_{\bold{g},D}$ is a direct consequence of (5.18). We consider $\langle S_{\bold{g},D}\bold{f},\bold{k}\rangle$ and observe that this 
is the left hand side of (5.18). Therefore (5.18) yields that
$$
\langle S_{\bold{g},D}\bold{f},\bold{k}\rangle=\big\langle\frac{1}{vol(D)}\sum_{h^!\in D^!}(V_{\bold{g}}\bold{g}\cdot V_{\bold{f}}\bold{k})(h^!)\big\rangle.
\eqno(5.19)
$$
and  we can write the right hand side of (5.19) in the following way:
$$
\sum_{d^!\in D^!}(V_{\bold{g}}\bold{g}\cdot V_{\bold{f}}\bold{h})(d^!)=\big\langle\sum_{d^!\in D^!} \langle\bold{g},\pi(d^!)\bold{g}\rangle\pi(d^!)\bold{f},\bold{h}\big\rangle.
\eqno(5.20)
$$ 
From this we derive the {\it Janssen representation} of $S_{\bold{g},D}$:
\medskip
{\bf 5.1.2. Proposition.} { \it Let $\bold{f},\bold{g},\bold{k}$ be in $M_1(\bold{R}^N)$ or in $\Cal{S}(\bold{R}^N)$. Then}
$$
S_{\bold{g},D}\bold{f}=\frac{1}{vol(D)}\sum_{d^!\in D^!}\langle\bold{g},\pi(d^!)\bold{g}\rangle\pi(d^!)\bold{f}. 
\eqno(5.21)
$$

The Janssen representation lies at the heart of the {\it duality theory} for Gabor frames, which was developed independently by three groups of researchers  [DaLaLa], [Ja4], [RoSh] at the same 
time. Later Feichtinger and Gr\"ochenig extended their $L_2(\bold{R})$ results. These results are the cornerstones of Gabor analysis. One important consequence of (5.21) is that the 
invertibility of $S_{\bold{g},D}$ becomes a question about the invertibility of an absolutely convergent sequence of time--frequency shifts:
$$
S_{\bold{g},D}=\frac{1}{vol(D)}\sum_{h^!\in D^!}\langle\bold{g},\pi(h^!)\bold{g}\rangle\pi(h^!). 
\eqno(5.22)
$$

In [GrLe] Gr\"ochenig and Leinert were able to link this fact with the spectral invariance of subalgebras of the quantum torus $C^*(D^!,\alpha)$. Later we observed that these kind of results are of great relevance for noncommutative geometry, especially some theorems of Connes on smooth quantum tori ([Lu2]).  

\smallskip
Observe that all our statements involving $M_1(\bold{R}^N)$ in this section also hold for $M_1^s(\bold{R}^N)$. Therefore by the time--frequency description of the smooth vectors of 
the vector Heisenberg group, i.e. the Schwartz class $\Cal{S}(\bold{R}^N)=\bigcap_{s\ge 0}M_1^s(\bold{R}^N)$, the spaces $\{M_1^s(\bold{R}^N)\,|\,s\ge 0\}$ provide a scale of Banach 
spaces that interpolates between the integrable vectors and the smooth vectors. 
 
\medskip
{\bf 5.2. Projective modules over quantum tori in the setting of Gabor analysis.}  In our treatment of the Janssen representation we emphasized the relevance of FIGA and that 
it already appeared in Rieffel's construction of projective modules over noncommutative tori. Now we want to exploit this link between projective modules over quantum tori 
and Gabor frames in more detail, see [Lu1], [Lu3], [Lu4] for further aspects of this topic.

\smallskip

The link between abstract quantum tori from Section 4 and Gabor analysis is provided by the choice of a particular representation for the quantum torus. Namely 
the operators $\{\pi(d):d\in D\}$ extend to a faithful involutive representation of the quantum torus $C^*(D,\alpha)$ on $L_2(\bold{R}^N)$. Consequentely the smooth subalgebras 
$C_1^s(D,\alpha)$ and $C^\infty(D,\alpha)$ become the following classes of operator algebras:

$$
C_1^s(D,\alpha)=\{\sum_{h\in D}a_h\pi(h)\,|\,\sum_{h\in D}|a_h|(1+|h|^2)^{s/2}<\infty\}
\eqno(5.23)
$$
and $C^\infty(D,\alpha)=\bigcap_{s\ge 0}C_1^s(D,\alpha)$. By the integrated Heisenberg commutation relations the 2-cocycle $\alpha$  in this representation take the following 
form: $\alpha(h,k)=e^{2\pi i h_\omega\cdot k_x}$ for $h=(h_x,h_\omega)$ and $k=(k_x,k_\omega)$. 

\smallskip

In other words $C_1^s(D,\alpha)$ consists of the image of the integrated representation of $\ell_1^s(D)$ with respect to the faithful involutive representation generated by 
the time--frequency shifts $\pi(h)$. We denote the integrated representation of $\bold{a}\in \ell_1^s(D)$ by $\pi_D(\bold{a})$:
$$
\pi_D(\bold{a})=\sum_{h\in D}a_h\pi(h).
\eqno(5.24)
$$

In this setting the left action of $C_1^s(D,\alpha)$ on $\bold{g}\in M_1^s(\bold{R}^N)$ becomes

$$
\bold{a}\cdot \bold{g}=\pi_D(\bold{a})\cdot\bold{g}= \sum_{h\in D}a_h\pi(h)\bold{g},
\eqno(5.25)
$$
for $\bold{a}\in\ell_1^s(D,\alpha)$ and similarly for $C^\infty(D,\alpha)$ acting on the left on $\Cal{S}(\bold{R}^N)$ for $\bold{a}\in\Cal{S}(D,\alpha)$. 

\smallskip

The algebra-valued product on $C_1^s(D,\alpha)$ turns out to be 
$$
{_D}\langle\bold{f},\bold{g}\rangle=\sum_{h\in D}\langle \bold{f},\pi(h)\bold{g}\rangle\pi(h),
\eqno(5.26)
$$
for $\bold{f},\bold{g}\in M_1^s(\bold{R}^N)$ and similarly on $C^\infty(D,\alpha)$ for $\bold{f},\bold{g}\in \Cal{S}(\bold{R}^N)$. 

\smallskip

Note that the left action (5.25) is just the synthesis operator of the Gabor system $\bold{G}(\bold{g},D)$ for $\bold{g}\in M_1^s(\bold{R}^N)$ or in 
$\Cal{S}(\bold{R}^N)$, and that the scalar product is the integrated representation of the sequence obtained from the analysis mapping of the Gabor system $\bold{G}(\bold{g},D)$. 

\medskip

{\bf 5.1.3. Proposition.} {\it Let $D$ be an embedded lattice of $\bold{R}^{2N}$. Then $M_1^s(\bold{R}^N)$ and $\Cal{S}(\bold{R}^N)$ become a
left Hilbert $C^*(D,\alpha)$--module $_DV$ after completing by the norm $_{D}\|\bold{f}\|=\|{_D}\langle\bold{f},\bold{f}\rangle\|^{1/2}$.} 

\medskip

We refer the reader to [Ri5], [Lu3], [Lu4] for a proof and generalizations of the last propostion. The most important operators on  Hilbert $C^*$-modules are 
the {\it rank--one Hilbert module operators}, which in our case are defined by 

$$
\Theta^D_{\bold{g},\bold{k}}\bold{f}:={_D}\langle\bold{f},\bold{g}\rangle\cdot\bold{k}=\sum_{h\in D}\langle \bold{f},\pi(h)\bold{g}\rangle\pi(h)\bold{k}.
\eqno(5.27)
$$

Operators of the form $\Theta^D_{\bold{g},\bold{k}}\bold{f}$ are called {\it Gabor frame--type operators} for a given Gabor system $\bold{G}(\bold{g},D)$ and are denoted by 
$S_{\bold{g},\bold{k},D}$. These operators appear naturally in the discrete reconstruction formulas. 

\smallskip
In the following we want to discuss the Rieffel-Morita equivalence of $C^*(D,\alpha)$ and $C^*(D^!,\overline{\alpha}^!)$. Recall that two $C^*$-algebras $\Cal{A}$ and $\Cal{B}$ are called {\it Rieffel-Morita equivalent}, if there exists an $\Cal{A}$-$\Cal{B}$-{\it equivalence bimodule} ${_\Cal{A}}V{_\Cal{B}}$ such that:
 \smallskip
   (i) ${_\Cal{A}}V{_\Cal{B}}$ is a full left Hilbert $\Cal{A}$-module and a full right Hilbert $\Cal{B}$-module;
 \smallskip  
   (ii) for all $f,g\in {_\Cal{A}}V{_\Cal{B}}, A\in\Cal{A}$ and $B\in\Cal{B}$ we have that
     $$
       \langle A\cdot f,g\rangle{_\Cal{B}}=\langle f,A^*\cdot g\rangle{_\Cal{B}}~~\text{and}~~{_\Cal{A}}\langle f\cdot B,g\rangle={_\Cal{A}}\langle f,g\cdot B^*\rangle;
     $$
  \smallskip  
    (iii) for all $f,g,h\in {_\Cal{A}}V{_\Cal{B}}$,
       $$
       {_\Cal{A}}\langle f,g\rangle\cdot h=f\cdot\langle g,h\rangle{_\Cal{B}}.
     $$
We refer the reader to [RaWi] for an extensive discussion of Rieffel-Morita equivalence. 
\smallskip
We continue with a discussion of the right module structure on $M_1^s(\bold{R}^N)$ or $\Cal{S}(\bold{R}^N)$. Since the quantum torus is only Morita equivalent to the opposite algebra torus
 $C^*(D^!,\overline{\alpha}^!)$, the action and the scalar product $\langle.,.\rangle_{D^!}$ differ from those in the case of $C^*(D,\alpha)$. More precisely, for $\bold{b}\in\ell^s_1(D^!)$, 
$\bold{f},\bold{g}$ in $M_1^s(\bold{R}^N)$ we define a {\it right action} of $C_1^s(D^!,\overline{\alpha}^!)$ on $M_1^s(\bold{R}^N)$ by

$$
\bold{g}\cdot\bold{b}:=\bold{g}\cdot\bold{\pi}_{D^!}(\bold{b})=\frac{1}{vol(D)}\sum_{h^!\in D^!}\pi(h^!)\bold{f}\overline{{b}(h^!)},
\eqno(5.28)
$$
and the $C_1^s(D^!,\overline{\alpha}^!)$-valued scalar product by 
$$
\langle\bold{f},\bold{g}\rangle_{D^!}=\sum_{h^!\in D^!}\pi(h^!)^*\langle\pi(h^!)\bold{g},\bold{f}\rangle.
\eqno(5.29)
$$

As indicated above this yields also a right action and a $C^\infty(D^!,\overline{\alpha}^!)$-valued scalar product in the case of $\Cal{S}(\bold{R}^N)$. 
Consequentely, we get a right Hilbert $C^*(D^!,\overline{\alpha}^!)$-module $V_{D^!}$ after completing $M_1^s(\bold{R}^N)$ or $\Cal{S}(\bold{R}^N)$ by 
$\|\bold{f}\|_{D^!}=\|\langle\bold{f},\bold{f}\rangle_{D^!}\|^{1/2}$. 

\smallskip
 
The main result of \cite{Ri5} asserts that these two Hilbert $C^*$-modules form actually an {\it equivalence bimodule}, 
and the two scalar products 
$_{D^!}\langle.,.\rangle$ and $\langle.,.\rangle_{D^!}$ satisfy the {\it associativity equation}:

$$
{_D}\langle\bold{f},\bold{g}\rangle\cdot\bold{k}=\bold{f}\cdot\langle\bold{g},\bold{k}\rangle_{D^!}
\eqno(5.30)
$$
for $\bold{f},\bold{g},\bold{k}$ in $M_1^s(\bold{R}^N)$ or $\Cal{S}(\bold{R}^N)$, respectively. Recall our discussion of the Janssen representation of 
$S_{\bold{g},D}$. If you consider the Janssen representation of a Gabor frame--type operator $S_{\bold{g},\bold{k},D}$, then you get the associativity condition (5.30). 

\medskip

{\bf 5.1.4. Proposition.} {\it Let $D$ be an embedded lattice of $\bold{R}^{2N}$. $M_1^s(\bold{R}^N)$ and $\Cal{S}(\bold{R}^N)$ become 
an $C^*(D,\alpha)$--\,$C^*(D^!,\overline{\alpha}^!)$ equivalence bimodule $_DV_{D^!}$after completing by the norm $_{D}\|\bold{f}\|=\|{_D}\langle\bold{f},\bold{f}\rangle\|^{1/2}$.}

\medskip 

The involutive Banach algebras $C_1^s(D,\alpha)$ and the smooth noncommutative torus $C^\infty(D,\alpha)$ are spectrally invariant subalgebras of the quantum torus 
$C^*(D,\alpha)$, i.e. if $A\in C_1^s(D,\alpha)$ or $C^\infty(D,\alpha)$ is invertible in $C^*(D,\alpha)$, then $A^{-1}$ is in $C_1^s(D,\alpha)$ or $C^\infty(D,\alpha)$, 
respectively. The spectral invariance of $C_1^s(D,\alpha)$ was recently proved by Gr\"ochenig and Leinert in [GrLe]. For $C^\infty(D,\alpha)$ this is a famous theorem of Connes [Co1]. 

\smallskip

In [Co2] it is demonstrated that if there exists a subspace $V_0$ of the equivalence bimodule $_DV_{D^!}$, that is invariant under the left and 
right actions of $\Cal{A}_0$ and $\Cal{B}_0$, and such that the scalar products evaluated for elements in $V_0$ are elements of spectrally invariant subalgebras $\Cal{A}_0$ and $\Cal{B}_0$ of $C^*(D,\alpha)$ and $C^*(D^!,\overline{\alpha}^!)$, then 
$V_0$ is a projective left $\Cal{A}_0$-module and projective right $\Cal{B}_0$-module respectively. We call the two involutive algebras $\Cal{A}_0$ and $\Cal{B}_0$ {\it Rieffel--Morita equivalent}, if there exists such an equivalence bimodule $V_0$. 

\smallskip

The spaces $M_1^s(\bold{R}^N)$ and the noncommutative Wiener algebras $C_1^s(D,\alpha)$ fulfill these requirements. Therefore we get 

\medskip

{\bf 5.1.5. Theorem.} {\it The noncommutative Wiener algebras $ C_1^s(D,\alpha)$ and $ C_1^s(D^!,\overline{\alpha}^!)$ are Rieffel-Morita equivalent through the modulation space $M_1^s(\bold{R}^N)$ for all $s\ge 0$. 
Consequently,  the noncommutative smooth tori $ C^\infty(D,\alpha)$ and $ C^\infty(D^!,\overline{\alpha}^!)$ are Rieffel-Morita equivalent through the Schwartz space $\Cal{S}(\bold{R}^N)$.}

\smallskip

In [Lu4] we generalize these results to other classes of spectrally invariant subalgebras of $C^*(D,\alpha)$ with subexponential growth. 
\bigskip


\centerline{\bf \S 6. Quantum theta functions and their functional equations.}
\medskip

For the treatment of quantum theta functions we have to introduce the class of generalized Gaussians, which appear prominently in various areas, e.g. in quantum optics as squeezed states, in harmonic analysis [Fo], the work of Littlejohn on semi-classical mechanics in [Li], and in de Gosson's work on symplectic capacity as measure of uncertainty in quantum mechanics [Go].

\smallskip

In the description of the $II_T$--model  in 2.2 above, we have implicitly used
the {\it Siegel upper--half plane} $\Cal{S}_N$ which is the space of all matrices of the form $T=\roman{Re}\,T+i\roman{Im}\,T$, where $\roman{Re}\,T,\roman{ImT}$ are real symmetric $N\times N$-matrices, and $\roman{Im}\,T$ is positive definite. Then we define a {\it generalized Gaussian} by
$$
  \bold{g}_{T}(x)=e^{-\langle Tx,x\rangle}=e^{-\langle(\roman{Re}\,T+i\roman{Im}\,T)x,x\rangle}~~\text{for}~~x\in\bold{R}^N.
\eqno(6.1)  
$$
In his work on quantum theta functions Manin [Ma2,Ma3,Ma4] has calculated  the Wigner transform of generalized Gaussians, which is a well--known result in quantum optics and quantum mechanics. 

\smallskip

Recall that the {\it Wigner distribution} $W(\bold{g},\bold{g})$ of a function $g$ is given by 
$$
  W(\bold{g},\bold{g})(x,\omega)=\int_{\bold{R}^N} \bold{g}(x+\tfrac{t}{2})\overline{\bold{g}}(x-\tfrac{t}{2})e^{-2\pi it\cdot\omega}dt
\eqno(6.2)
$$
and that the Wigner distribution is just a short--time Fourier transform in disguise:
$$
   W(\bold{g},\bold{g})(x,\omega)=2^Ne^{4\pi ix\cdot\omega}V_{\tilde \bold{g}}\bold{g}(2x,2\omega),
\eqno(6.3)   
$$
where $\tilde \bold{g}(x)=\bold{g}(-x)$ denotes the reflection of $\bold{g}$ with respect to the origin, see eq. (3.8) and (3.9).
\par
Probably, one of the earliest calculations of the Wigner transform of a generalized Gaussian has been published by the electrical engineer Bastiaans [Ba]. 

\smallskip

{\bf 6.1. Lemma.} {\it Let $\bold{g}_{T}$ be the generalized Gaussian for a $T\in\Cal{S}_N$. Then the Wigner transform of $\bold{g}_{T}$ is given by} 
$$
  W(\bold{g}_{T},\bold{g}_{T})(z)=(\roman{det}\,T)^{-1/2}e^{-H_T(z,z)}
\eqno(6.4)  
$$
{\it where $H_T(z,z)=G_Tz\cdot z$, and }                    
$$
  G_T=
\left(
\matrix 
 \roman{Re}\,T+ \roman{Im}\,T( \roman{Re}\,T)^{-1}\roman{Im}\,T & \roman{Im}\,T(\roman{Re}\,T)^{-1}\\
(\roman{Re}\,T)^{-1}\roman{Im}\,T & (\roman{Re}\,T)^{-1}
\endmatrix \right)
\eqno(6.4)
$$
{\it The matrix $G_T$ is symplectic and has the following factorization: $G_T=S^TS$ with}
$$
S=
\left(
\matrix 
(\roman{Re}\,T)^{1/2} & 0\\
(\roman{Re}\,T)^{-1/2}\roman{Im}\,T & (\roman{Re}\,T)^{-1/2}
\endmatrix \right).
\eqno(6.5)
$$
\smallskip

A generalized Gaussian $\bold{g}_T$ is an element of $\Cal{S}(\bold{R}^N)$, because it belongs also to $M^1_s(\bold{R}^N)$ for any $s\ge 0$. 
The last observation follows from a well--known characterization of $M^1_s(\bold{R}^N)$ in terms of the Wigner transform. Namely, 
$\bold{g}\in M^1_s(\bold{R}^N)$ if and only if $W(\bold{g},\bold{g})$ is in $L^1_s(\bold{R}^N)$. An elementary computation yields 
that the Wigner distribution of $\bold{g}_T$ is in $L^1_s(\bold{R}^N)$ for any $s\ge 0$. We combine this observation and its interesting 
consequence in the next lemma.

\smallskip

{\bf 6.2. Lemma.} {\it Let $D$ be an embedded sublattice of $\bold{R}^{2N}$. Since $\bold{g}_T$ is in $M^1_s(\bold{R}^N)$ for all $s$, 
 we have that $(\langle \bold{g}_T,\pi(h)\bold{g}_T\rangle)_{h\in D}$ is in $\ell^1_s(D)$ and ${_D}\langle\bold{g}_T,\bold{g}_T\rangle$
defines an element of $C_1^s(D,\alpha)$.}

\smallskip

A natural interpretation of the quantum theta functions studied in [Ma2,Ma3,Ma4] comes from the study of the Gabor system $\bold{G}(\bold{g}_T,D)$ for 
a general Gaussian $\bold{g}_T$ and an embedded lattice $D$ of $\bold{R}^{2N}$. The Gabor frame operator $S_{\bold{g}_T,D}$ looks as follows 
$$
S_{\bold{g}_T,D}\bold{f}=\sum_{h\in D}\langle\bold{f},\pi(h)\bold{g}_T\rangle\pi(h)\bold{g}_T,
\eqno(6.6)
$$
for $\bold{f}$ in $L_2(\bold{R}^N)$. If we consider the Gabor frame operator on $M_1^s(\bold{R}^N)$, then the Janssen representation of $S_{\bold{g}_T,D}$ exists and turns out to be a quantum theta function:
$$
S_{\bold{g}_T,D}=\frac{1}{vol(D)}\sum_{h^!\in D^!}\langle\bold{g}_T,\pi(h^!)\bold{g}_T\rangle\pi(h^!). 
\eqno(6.7)
$$
By Lemma 6.1 the Janssen representation can be rewritten as
$$
S_{\bold{g}_T,D}=\frac{1}{vol(D)}\sum_{h^!\in D^!}e^{-\pi H_T(h^!,h^!)}\pi(h^!).
\eqno(6.8)
$$ 
Therefore, the superposition of time--frequency shifts in the preceding equation  
$$
\sum_{h^!\in D^!}e^{-\pi H_T(h^!,h^!)}\pi(h^!)
\eqno(6.9)
$$
is an element of $C^s_1(D^!,\overline{\alpha}^!)$ for each $s$ and consequently of $C^\infty(D^!,\overline{\alpha}^!)$. In [Ma5] the operator of (6.9) was denoted by $\Theta_{D^!}$ 
and it was noted that $\Theta_{D^!}$ is a {\it quantum theta function} in $C^s_1(D^!,\alpha)$ and consequently in $C^\infty(D^!,\overline{\alpha}^!)$.  More explicitely, $\Theta_{D^!}$ satisfies the following functional 
equations:
$$
\forall\, h^!\in D^!,\ C_{h^!} \pi_{h^!}x_{h^!}^*(\Theta_{D^!})=\Theta_{D^!},
\eqno(6.10)
$$
where
$$
C_{h^!}=e^{-\frac{\pi}{2}\,H(\underline{h^!},\underline{h^!})},\ 
x_{h^!}^*(\pi(h^!)) = e^{-\pi\,H(\underline{h^!},\underline{h^!})}\pi(h^!).
\eqno(6.11)
$$
In the Model $II_T$ the quantum theta function $\Theta_{D^!}$ becomes
$$
\Theta_D\,\bold{1}(x)=
\sum_{h^!\in D} 
e^{-\pi H(\underline{h^!},\underline{h^!})-\pi H(\underline{x},\underline{h^!})} 
\eqno(6.12)
$$ 
where $\bold{1}$ is the vacuum vector in the model II${}_T$ represented by
the function identically equal to 1.

\smallskip

The function  $\Theta_{D^!} \bold{1}$ is {\it complex conjugate} to the classical theta function
corresponding to a principal polarization of the complex torus $\bold{C}^g/D^!$. Notice that this complex torus
is embedded into (the space of points of) the algebraic torus $T(D^!,1)(\bold{C}) = \roman{Hom}\,(D^!,\bold{C}^*)$
as its compact subtorus $\roman{Hom}\,(D^!,\bold{C}_1^*)$, see [Mu]. 

\smallskip

Moreover, the functional equation for quantum thetas in [Ma5] is an expression for the Janssen representation of 
$S_{\bold{g}_T,D}\bold{f}$ for $\bold{f}=\bold{g}_T$:

$$
\sum_{h\in D}
e^{-\pi H(\underline{h},\underline{h})-\pi H(\underline{x},\underline{h})}=
\frac{1}{vol(D)}\sum_{h^!\in D^!}
e^{-\pi H(\underline{h^!},\underline{h^!})-\pi H(\underline{x},\underline{h^!})}.
\eqno(6.13)
$$
as functions of $x\in \bold{R}^{2N}.$

\smallskip
Since the Janssen representation is a consequence of FIGA one may also consider FIGA for $\bold{f}_1=\bold{f}_2=\bold{g}_1=\bold{g}_2=\bold{g}_T$ 
as a functional equation for quantum thetas. In the case when $N=1$ and $g_T$ is the standard Gaussian, Schempp already noted in [Sc] that this kind 
of identities are linked with theta functions.
\smallskip
Note that the quantum thetas $\Theta_D$ in [Ma5] are defined by ${_D}\langle \bold{f}_T,\bold{f}_T\rangle$, where $\bold{f}_T$ 
is our $\bold{g}_{iT}$.

\smallskip

We close this section with a few words on the case $D=a\bold{Z}\times b\bold{Z}$. In this case $H_T(z,w)=z\overline{w}$ for $w,z\in\bold{C}$. 
Consequently $\bold{g}_T$ is just the standard Gaussian $\bold{g}_0(x)=e^{-\pi x^2}$. Therefore the quantum thetas $\Theta_{a\bold{Z}\times b\bold{Z}}$ 
are related to the Gabor systems $\bold{G}(\bold{g}_0,a\bold{Z}\times b\bold{Z})$. A deep result of Lyubarski and Seip obtained independely in [Ly], [Se] says that $\bold{G}(\bold{g}_0,a\bold{Z}\times b\bold{Z})$ is a Gabor frame for $L_2(\bold{R}^N)$ if and only if $ab<1$. An important result of Feichtinger and Gr\"ochenig [FeGr1] asserts that Gabor frames $\bold{G}(g,D)$ for $L_2(\bold{R}^N)$ with $\bold{g}$ in $M_1^s(\bold{R}^N)$ or $\Cal{S}(\bold{R}^N)$ are {\it Banach frames} for the class of modulation spaces $M_{p,q}^m(\bold{R}^N)$. Therefore $\bold{G}(\bold{g}_0,a\bold{Z}\times b\bold{Z})$ is a  tight Banach frame for $M_{p,q}^m(\bold{R}^N)$ if and only if $ab<1$, i.e. 
$$
\|\bold{f}\|_{M_{p,q}^m}=\big(\sum_{l\in\bold{Z}}\big(\sum_{k\in\bold{Z}}|\langle\bold{f},\pi(ak,bl)\bold{g}_0\rangle|^p\big)^{q/p}\big)^{1/q}<\infty
\eqno(6.14)
$$
for all $\bold{f}$ in $M_{p,q}^m(\bold{R}^N)$ and $p,q\in[0,\infty]$. We want to emphasize that the preceding equation provides a description 
in terms of Gabor coefficients of ``good" Gabor frames.

\smallskip

Observe that $\bold{G}(\bold{g}_0,a\bold{Z}\times b\bold{Z})$ is a Gabor frame for $L_2(\bold{R})$ if and only if the Janssen representation of the 
Gabor frame operator is invertible. Since the Janssen representation of $S_{\bold{g}_0,a\bold{Z}\times b\bold{Z}}$ is the quantum theta 
$\Theta_{\tfrac{1}{b}\bold{Z}\times\tfrac{1}{a}\bold{Z}}$, we are thus able to 
produce a precise criterium for the invertiblity of the quantum thetas 
$\Theta_{\tfrac{1}{b}\bold{Z}\times\tfrac{1}{a}\bold{Z}}$. Moreover the spectral invariance of $C_1^s(\tfrac{1}{b}\bold{Z}\times\tfrac{1}{a}\bold{Z},\alpha)$ in $C^*(\tfrac{1}{b}\bold{Z}\times\tfrac{1}{a}\bold{Z},\alpha)$ implies that $\Theta_{\tfrac{1}{b}\bold{Z}\times\tfrac{1}{a}\bold{Z}}^{-1}$ is in $C_1^s(\tfrac{1}{b}\bold{Z}\times\tfrac{1}{a}\bold{Z},\alpha)$ if and only if $ab<1$. These observations provide a new approach to the projections in [Bo] and in addition clarifies the connection between quantum thetas and 
these projections in quantum tori. 

\smallskip

{\bf 6.3. Proposition.} {\it The quantum theta $\Theta_{\tfrac{1}{b}\bold{Z}\times\tfrac{1}{a}\bold{Z}}$ is invertible if and only if $ab<1$.}

\smallskip

The construction of projections in higher dimensional quantum tori using Gabor analysis will be addressed by one of us in a subsequent publication. The Gabor systems $\bold{G}(\bold{g}_T,D)$ will also be studied in more detail.

\bigskip

{\bf Acknowledgement.} Large parts of this manuscript were written during a stay of F.~L. at the 
Max Planck Insitute for Mathematics at Bonn. F.~L. would like to express his gratitude for hospitality and 
excellent working conditions. In addition F.~L. acknowledges the support from the EU-project MEXT-CT-2004-517154 and the Marie Curie Outgoing Fellowship PIOF 220464.

\newpage

\centerline{\bf Bibliography}

\medskip

[Ba] M.~J. {B}astiaans.
{\it {W}igner distribution function and its application to first-order
optics.}
{{J}. {O}pt. {S}oc. {A}mer.}, 69(12) (1979), 1710--1716. 
\smallskip
[Bo1] F.~Boca. {\it Projections in rotation algebras and theta functions.}
Comm. Math. Phys., 202 (1999), 325--357.

\smallskip

[Ca] P.~{C}artier. {\it
 {Q}uantum mechanical commutation relations and theta functions.}
In: {A}lgebraic {G}roups and {D}iscontinuous {S}ubgroups ({P}roc.
  {S}ympos. {P}ure {M}ath., {B}oulder, {C}olo., 1965), {A}mer.
  {M}ath. {S}oc., {P}rovidence, {R}.{I}., 1966,   361--383.

\smallskip

[ChKi] E.~{C}hang {Y}oung and H.~{K}im.
{\it {T}heta vectors and quantum theta functions.} 
{{J}. {P}hys. {A}}, 38(19) (2005), 4255--4263.

\smallskip

[Co1] A.~{C}onnes. {\it $C^*$-alg\`ebres et g\'eom\'etrie diff\'erentielle.}
{{C}. {R}. {A}cad. {S}ci. {P}aris {S}\'er. {A}-{B}},  290(13) (1980), {A}599--{A}604.

\smallskip

[Co2] A.~{C}onnes. {\it An analogue of the {T}hom isomorphism for crossed products of a {$C^{*} $}-algebra by an action of {$\bold{R}$}}
{Adv. in Math.},  39(1) (1981), 31--55.

\smallskip

[Co3] A.~Connes. {\it Noncommutative geometry.} Academic Press,
1994.

\smallskip

[DaGrMe] I.~{D}aubechies, A.~{G}rossmann, and Y.~{M}eyer.
{\it Painless nonorthogonal expansions. }
 {J}. {M}ath. {P}hys., 27(5) (1986), 1271--1283.

\smallskip

[Da] I.~{D}aubechies.
{\it The wavelet transform, time--frequency localization and signal analysis.}
{I}{E}{E}{E} {T}rans. {I}nf. {T}heory, 36(5) (1990), 961--1005.

\smallskip

[DaLaLa] I.~{D}aubechies, H.~J. {L}andau, and Z.~{L}andau.
{\it {G}abor time--frequency lattices and the {W}exler--{R}az identity.}
{{J}. {F}ourier {A}nal. {A}ppl.}, 1(4) (1990), 437--478.

\smallskip

[DuSc] R.~J.~Duffin, A.~C.~ Schaeffer, A. C.
{\it A class of nonharmonic {F}ourier series},
Trans. Amer. Math. Soc., 72 (1952), 341--366.

\smallskip

[EeK]  Ee C.--Y., H.~Kim. {\it Quantum thetas on noncommutative $T^d$ with general
embeddings.}  e--print math-ph/0709.2483

\smallskip

[FaSc]  {D}.~{C}.~{F}arden and L.~L.~ {S}charf.
{\it {A} unified framework for the {S}ussman, {M}oyal, and {J}anssen formulas.}
{I}{E}{E}{E} {S}ignal {P}rocessing {M}agazine, 124 (2006), 124--125.

\smallskip

[Fe1] H.~G. {F}eichtinger. {\it {O}n a new {S}egal algebra.}
{{M}onatsh. {M}ath.}, 92 (1981), 269--289.

\smallskip

[Fe2] H.~G. {F}eichtinger. {\it {M}odulation spaces on locally compact {A}belian groups.}
Technical report, {J}anuary 1983.

\smallskip

[Fe3] H.~G. {F}eichtinger. {\it {M}odulation {S}paces: {L}ooking {B}ack and {A}head.}
{S}ampl. {T}heory {S}ignal {I}mage {P}rocess., 5(2) (2006), 109--140.

\smallskip

[FeGr] H.~G. {F}eichtinger and K.~{G}r{\"o}chenig. {\it {B}anach spaces related to integrable group representations and their
atomic decompositions, {I}.}
{{J}. {F}unct. {A}nal.}, 86 (1989), 307--340.

\smallskip

[FeGr1] H.~G. {F}eichtinger and K.~{G}r{\"o}chenig. {\it {G}abor frames and time--frequency analysis of distributions.}
{{J}. {F}unct. {A}nal.}, 146(2) (1997), 464--495 .

\smallskip

[FeKo] H.~G. {F}eichtinger and W.~{K}ozek.
{\it {Q}uantization of {T}{F} lattice-invariant operators on elementary {L}{C}{A} groups.}
In H.~{F}eichtinger and T.~{S}trohmer, editors, {{G}abor {A}nalysis and {A}lgorithms. {T}heory and {A}pplications.}, 
{A}pplied and {N}umerical {H}armonic {A}nalysis, pages 233--266, 452--488, {B}oston, {M}{A}, 1998. {B}irkh{\"a}user {B}oston.

\smallskip

[FeLu] H.~G. {F}eichtinger and F.~{L}uef.
{\it {W}iener amalgam spaces for the fundamental identity of {G}abor analysis.}
{{C}ollect. {M}ath.}, 57 (2006), 233--253.
e--Print math.FA/0503364

\smallskip

[Fo] G.~B. {F}olland.
{\it {H}armonic {a}nalysis in {p}hase {s}pace}.
{P}rinceton {U}niversity {P}ress, {P}rinceton, {N}.{J}., 1989.

\smallskip

[FrLa] M.~{F}rank and D.~R. {L}arson.
{\it {F}rames in {H}ilbert ${C}^*$-modules and ${C}^*$-algebras.} 
{{J}. {O}perator {T}heory}, 48 (2002), 273--314.

\smallskip

[Ga] D.~{G}abor.
{\it {T}heory of communication.}
{{J}. {I}{E}{E}}, 93(26) (1946), 429--457.

\smallskip

[Go] M.~de~{G}osson.
{\it {S}ymplectic {G}eometry and {Q}uantum {M}echanics}, volume 166
  of { {O}perator {T}heory: {A}dvances and {A}pplications. {A}dvances in
  {P}artial {D}ifferential {E}quations.}, {B}irkh{\"a}user, 2006.

\smallskip

[Gr1] K.~{G}r{\"o}chenig.
{\it {A}spects of {G}abor analysis on locally compact abelian groups.}
In H.~{F}eichtinger and T.~{S}trohmer, editors, { {G}abor analysis
  and algorithms: {T}heory and {A}pplications}, pages 211--231.
  {B}irkh{\"a}user {B}oston, {B}oston, {M}{A}, 1998.
  
\smallskip

[Gr2] K.~{G}r{\"o}chenig.
{\it {F}oundations of {t}ime--{f}requency {a}nalysis}.
{A}ppl. {N}umer. {H}armon. {A}nal. {B}irkh{\"a}user {B}oston, {B}oston, {M}{A}, 2001.

\smallskip

[GrLe] K.~{G}r{\"o}chenig and M.~{L}einert.
{\it {W}iener's {l}emma for {t}wisted {c}onvolution and {G}abor {f}rames.}
{{J}. {A}mer. {M}ath. {S}oc.}, 17(1) (2004), 1--18.

\smallskip

[GrLu] K.~{G}r{\"o}chenig and F.~{L}uef.
{\it {T}he topological stable rank of projective modules over noncommutative tori.}
Preprint, 2008.

\smallskip

[Ho] R.~{H}owe.
{\it {O}n the role of the {H}eisenberg group in harmonic analysis.}
{B}ull. {A}m. {M}ath. {S}oc., 3(2) (1980), 821--843.

\smallskip
[Ig] J.-i. {I}gusa.
{\it {T}heta functions.}
{D}ie {G}rundlehren der mathematischen {W}issenschaften. {B}and 194.{S}pringer, 1972.

\smallskip

[Ja1] A.~J. E.~M. {J}anssen.
{\it {G}abor {R}epresentation of {G}eneralized {F}unctions.}
{{J}. {M}ath. {A}nal. {A}ppl.}, 83 (1981), 377--394.

\smallskip

[Ja2] A.~J. E.~M. {J}anssen.
{\it{W}eighted {W}igner distributions vanishing on lattices.}
{{J}. {M}ath. {A}nal. {A}ppl.}, 80 (1981), 156--167.

\smallskip

[Ja3] A.~J. E.~M. {J}anssen.
{\it {B}argmann transform, {Z}ak transform, and coherent states.}
{{J}. {M}ath. {P}hys.}, 23(5) (1982), 720--731.

\smallskip

[Ja4] A.~J. E.~M. {J}anssen.
{\it {D}uality and biorthogonality for {W}eyl-{H}eisenberg frames.}
{{J}. {F}ourier {A}nal. {A}ppl.}, 1(4) (1995), 403--436.

\smallskip

[Ko] M.~Kontsevich.  {\it Deformation quantization of algebraic varieties.}
e-Print math.AG/0106006

\smallskip

[Li] R.~G. {L}ittlejohn.
{\it {T}he semiclassical evolution of wave packets.}
{{P}hys. {R}ep.}, 138(4-5) (1986), 193--291.

\smallskip

[Lu1] F.~{L}uef.
{\it {G}abor {A}nalysis meets {N}oncommutative {G}eometry}.
PhD thesis, {U}niversity of {V}ienna, {N}ov. 2005.

\smallskip

[Lu2] F.~{L}uef.
{\it {O}n spectral invariance of non-commutative tori.}
In {{O}perator {T}heory, {O}perator {A}lgebras, and
  {A}pplications}, volume 414, {A}merican {M}athematical
  {S}ociety, (2006), 131--146.
e-Print math.OA/0603139 

\smallskip

[Lu3] F.~{L}uef.
{\it {G}abor analysis, noncommutative tori and {F}eichtinger's algebra.}
In {{G}abor and {W}avelet {F}rames}, volume~10 of {{I}{M}{S} {L}ecture {N}otes {S}eries}. {W}orld {S}ci.{P}ub., (2007), 77--106.
e--Print math.FA/0504146 

\smallskip

[Lu4] F.~{L}uef.
{\it {P}rojective modules over non-communtative tor are multi-window {G}abor frames for modulation spaces.}
e--Print math.OA/0807.3170.

\smallskip

[Ly] Y.~I. {L}yubarskij.
{\it {F}rames in the {B}argmann space of entire functions.} In {{E}ntire and subharmonic functions}, volume~11 of {{A}dv.
  {S}ov. {M}ath.}, pages 167--180. {A}merican {M}athematical {S}ociety
  ({A}{M}{S}), {P}rovidence, {R}{I}, 1992.
  
\smallskip

[Ma1] Yu.~Manin. {\it  Quantized theta--functions.} In: Common
Trends in Mathematics and Quantum Field Theories (Kyoto, 1990), 
Progress of Theor. Phys. Supplement, 102 (1990), 219--228.

\smallskip

[Ma2] Yu.~Manin. {\it Mirror symmetry and quantization of abelian varieties.}
In: Moduli of Abelian Varieties, ed. by C.~Faber et al.,
Progress in Math., vol. 195, Birkh\"auser, 2001, 231--254.
e--Print  math.AG/0005143

\smallskip

[Ma3] Yu.~Manin. {\it Theta functions, quantum tori and Heisenberg groups}.
Lett. in Math. Physics, 56:3 (2001), 295--320.
e--Print math.AG/001119

\smallskip

[Ma4] Yu.~Manin. {\it Real multiplication and noncommutative geometry (ein {A}lterstraum)}
{ In: The legacy of Niels Henrik Abel}, Springer,
Berlin, 2004, 685--727.
e--Print math.AG/0202109

\smallskip

[Ma5] Yu.~Manin. {\it Functional equations for quantum theta functions}.
Publ. Res. Inst. Math. Sci., 40:3 (2004), 605--624.
e--Print math.QA/0307393

\smallskip

[Mu] D.~Mumford (with M.~Nori and P.~Norman). 
{\it Tata Lectures on Theta III.} Progress in Math., vol.~97,
Birkh\"auser, 1991.

\smallskip

[Po] A.~{P}olishchuk.
{\it {A}belian varieties, theta functions and the {F}ourier
  transform}, volume 153 of { {C}ambridge {T}racts in {M}athematics}.
{C}ambridge {U}niversity {P}ress, {C}ambridge, 2003.

\smallskip 
[RaWi] I.~{R}aeburn and D.~P. {W}illiams.
{\it {M}orita {E}quivalence and {C}ontinuous-trace $C^*$-algebras.} 
{A}merican {M}athematical {S}ociety ({A}{M}{S}), {R}{I}, 1998. 

\smallskip

[Re1]
M. {R}eiter. {\it \"{U}ber den {S}atz von {W}eil--{C}artier.}
 {M}h. {M}ath., 86 (1978), 13--62.
 
 \smallskip
 
[Re2] M. {R}eiter.
{\it  {T}heta functions and symplectic groups.}
{M}h. {M}ath., 97 (1984), 219--232.
\smallskip

[Re3] H.~{R}eiter. {\it {M}etaplectic {G}roups and {S}egal {A}lgebras.} {S}pringer, {B}erlin, 1989.

\smallskip

[Ri1] M.~A. {R}ieffel.
{\it {C}*-algebras associated with irrational rotations.}
{{P}ac. {J}. {M}ath.}, 93 (1981) 415--429.

\smallskip

[Ri2] M.~A. {R}ieffel. {\it {T}he cancellation theorem for projective modules over irrational rotation ${C}^*$-algebras.}
{{P}roc. {L}ond. {M}ath. {S}oc., {I}{I}{I}. {S}er.}, 47 (1983) 285--302.

\smallskip

[Ri3] M.~A.~Rieffel. {\it Projective modules over higher--dimensional
non--commutative tori.} Can.~J.~Math., vol.~XL, No.~2 (1988), 257--338.

\smallskip

[Ri4] M.~A.~Rieffel. {\it Non--commutative tori --- a case
study of non--commutative differential manifolds.}
In: Cont.~Math., 105 (1990) 191--211.

\smallskip

[RiSch] M.~A.~Rieffel, A.~Schwarz. {\it Morita equivalence
of multidimensional non--commutative tori.} 
Int. J. Math. 10 (1999), 289--299. 
e--Print math.QA/9803057

\smallskip

[Ro] J.~Rosenberg. {\it Noncommutative variations on Laplace's equation.} 
e--Print math.OA/0802.4033, (2008).

\smallskip

[RoSh] A.~{R}on and Z.~{S}hen.
{\it {W}eyl-{H}eisenberg frames and {R}iesz bases in ${L}_2(\bold{R}^d)$}.
 {D}uke {M}ath. {J}., 89(2) (1997), 237--282.

\smallskip

[Sch2] A.~Schwarz. {\it Theta--functions on non--commutative tori.}
Lett. in Math. Physics, 58:1 (2001), 81--90. 
e--Print math.QA/0107186

\smallskip

[Sc] W.~{S}chempp.
{\it {R}adar ambiguity functions, the {H}eisenberg group, and holomorphictheta series.}
{{P}roc. {A}m. {M}ath. {S}oc.}, 92 (1984), 103--110.

\smallskip

[Se] K.~{S}eip.
{\it {D}ensity theorems for sampling and interpolation in the {B}argmann-{F}ock space. {I}.}
{{J}. {R}eine {A}ngew. {M}ath.}, 429 (1992), 91--106.

\smallskip

[Su] S.~M.~{S}ussman.
{\it {L}east square synthesis of radar ambiguity functions.}
{{I}{R}{E} {T}rans. {I}nform. {T}heory}, 8, (1962), 246--254.

\smallskip

[Vl] M.~{V}lasenko.
{\it {T}he graded ring of quantum theta functions for noncommutative torus
  with real multiplication.}, { {I}nt. {M}ath. {R}es. {N}ot.}, (2006), 1--19.

\smallskip

[We] A.~{W}eil.
{\it {S}ur certains groupes d'op{\'e}rateurs unitaires.}
{{A}cta {M}ath.}, 111 (1964), 143--211.

\bigskip

e-mail: franz.luef\@univie.ac.at
\smallskip
e-mail: manin\@mpim-bonn.mpg.de

\enddocument